\documentclass[color]{pjm}
\usepackage{color,bbold}
\usepackage{amsrefs}
\usepackage[nohams]{math}
\usepackage[arrow,curve,matrix]{xy}
\def\pair<#1>{{\ll}#1{\gg}}

\newcommand\gr{{\operatorname{\mathsf{growth}}}}
\newcommand\alt[1]{{\mathfrak A_{#1}}}
\newcommand\squish[1]{{\makebox[3em]{$\langle#1\rangle^\Gg$}}}
\newcommand\ssquish[1]{{\makebox[3em]{$\scriptstyle\langle#1\rangle^\Gg$}}}

\def\lar{\ar@{-}}
\def\tar{\ar@{->}} 
\def\aar{\ar@{.>}} 
\newcommand\rmI{{\mathrm{I}}}
\newcommand\rmII{{\mathrm{I\!I}}}
\newcommand\rmIII{{\mathrm{I\!I\!I}}}
\newcommand\6{\genfrac{}{}{0pt}1{\mathbb0}{\mathbb1}}
\newcommand\7{{\mathbb0}}
\newcommand\8{{\mathbb1}}
\newcommand\9{{\mathbb2}}
\newcommand\usp{{\mathbb p'}}

\CompileMatrices
\begin{document}
\title{Lie Algebras and Growth in Branch Groups}
\author{Laurent Bartholdi}
\date{April 18, 2002; Revised December 5, 2002 and \today}
\email{laurent@math.berkeley.edu}
\urladdr{\texttt{http://www.math.berkeley.edu/\char126laurent}}
\address{Dept of Mathematics, Evans Hall 970, U. C. Berkeley, CA
  94720-3840 U.S.A.}
\thanks{The author acknowledges support from the ``Swiss National Fund
  for Scientific Research'', and the Hebrew University of Jerusalem.}
\begin{abstract}
  We compute the structure of the Lie algebras associated to two
  examples of branch groups, and show that one has finite width while
  the other, the ``Gupta-Sidki group'', has unbounded width
  (Corollary~\ref{cor:gamma:rk}). This answers a question by Sidki.
  More precisely (Corollary~\ref{cor:gsgk}) the Lie algebra of the
  Gupta-Sidki group has Gelfand-Kirillov dimension
  $\log3/\log(1+\sqrt2)$.
  
  We then draw a general result relating the growth of a branch group,
  of its Lie algebra, of its graded group ring, and of a natural
  homogeneous space we call \emph{parabolic space}, namely the
  quotient of the group by the stabilizer of an infinite ray. The
  growth of the group is bounded from below by the growth of its
  graded group ring, which connects to the growth of the Lie algebra
  by a product-sum formula, and the growth of the parabolic space is
  bounded from below by the growth of the Lie algebra
  (see Theorem~\ref{thm:growth}).
  
  Finally we use this information to explicitly describe the normal
  subgroups of $\Gg$, the ``Grigorchuk group''.  All normal subgroups
  are characteristic, and the number $b_n$ of normal subgroups of
  $\Gg$ of index $2^n$ is odd and satisfies
  $\{\limsup,\liminf\}b_n/n^{\log_2(3)}=\{5^{\log_2(3)},\frac29\}$
  (see Corollary~\ref{cor:asympns}).
\end{abstract}
\maketitle

\section{Introduction}
The first purpose of this paper is to describe explicitly the Lie
algebra associated to the Gupta-Sidki group
$\GS$~\cite{gupta-s:burnside}, and show in this way that this group
is not of finite width (Corollary~\ref{cor:gamma:rk}). We shall
describe in Theorem~\ref{thm:gamma:struct} the Lie algebra as a graph,
somewhat similar to a Cayley graph, in a formalism close to that
introduced in~\cite{bartholdi-g:lie}.

We shall then consider another group, $\FG$, and show in
Corollary~\ref{cor:delta:rk} that although many similarities exist
between $\GS$ and $\FG$, the Lie algebra of $\FG$ does have
finite width.

These results follow from a description of group elements as ``branch
portraits'', exhibiting the relation between the group and its Lie
algebra. They lead to the notion of infinitely iterated ``wreath
algebras'', similar to wreath products of
groups~\cite{bartholdi:wreathalg}, to appear in a subsequent paper.

We shall show in Theorem~\ref{thm:growth} that, in the class of branch
groups, the growth of the homogeneous space $G/P$ (where $P$ is a
parabolic subgroup) is larger than the growth of the Lie algebra
$\Lie(G)$. This result parallels a lower bound on the growth of $G$ by
that of its graded group ring $\overline{\Bbbk G}$
(Proposition~\ref{prop:gpgr}).

Finally, we shall describe all the normal subgroups of the first
Grigorchuk group, using the same formalism as that used to describe
the lower central series. We confirm the description by Ceccherini et
al.\ of the low-index normal subgroups of
$\Gg$~\cite{ceccherini-s-t:grigns}.  It turns out that all non-trivial
normal subgroups are characteristic, and have finite index a power of
$2$.  Call $b_n$ the number of normal subgroups of index $2^n$
(Finite-index, non-necessarily-normal subgroups always have index a
power of $2$; this follows from $G$ being a $2$-torsion group.) Then
there are $3^k+2$ subgroups of index $2^{5\cdot2^k+1}$ and
$\frac293^k+1$ subgroups of index $2^{2^k+2}$; these two values are
extreme, in the sense that $b_n/n^{\log_2(3)}$ has lower limit
$5^{-\log_2(3)}$ and upper limit $\frac29$.  Also, $b_n$ is odd for
all $n$ (see Corollaries~\ref{cor:asympns} and~\ref{cor:allodd}).

\subsection{Philosophy}
One can hardly exaggerate the importance of Lie algebras in the study
of Lie groups. Lie subgroups correspond to subalgebras, normal
subgroups correspond to ideals; simplicity, nilpotence etc.\ match
perfectly. This is due to the existence of mutually-inverse functions
$\exp$ and $\log$ between a group and its algebra, and the
Campbell-Hausdorff formula expressing the group operation in terms of
the Lie bracket.

In the context of (discrete) $p$-groups and Lie algebras of
characteristic $p$, the correspondence is not so perfect. First, in
general, there is no exponential, and the best one can consider is the
degree-$1$ truncations
\[\exp(x) = 1+x+\mathcal O(x^2),\quad\log(1+x)=x+\mathcal O(x^2);\]
more terms would introduce denominators that in general are not
invertible; and no reasonable definition of convergence can be imposed
on $\F$. As a consequence, the group has to be subjected to a
filtration to yield a Lie algebra. Then there is no perfect bijection
between group and Lie-algebra objects.

However, the numerous results obtained in the area show that much can
be gained from consideration of these imperfect algebras. To name a
few, the theory of groups of finite width is closely related to the
classification of finite $p$-groups
(see~\cites{leedham-green:prop,shalev-z:finite-coclass}) and the
theory of pro-$p$-groups is intimately Lie-algebraic;
see~\cite{shalev:pro-p}, \cite{shalev:finite-p}*{\S8}
and~\cite{klass-lg-p:fw} with its bibliography. The solution to
Burnside's problems by Efim Zelmanov relies also on Lie algebras. The
results by Lev Kaloujnine on the $p$-Sylow subgroups of $\sym{p^n}$,
even if in principle independent, can be restated in terms of Lie
algebras in a very natural way~(see Theorem~\ref{thm:psylow}).

In this paper, I wish to argue that questions of growth, geometry and
normal subgroup structure are illuminated by Lie-algebraic
considerations.

\subsection{Notation}
We shall always write commutators as $[g,h]=g^{-1}h^{-1}gh$,
conjugates as $g^h=h^{-1}gh$, and the adjoint operators $\Ad(g)=[g,-]$
and $\ad(x)=[x,-]$ on the group and Lie algebra respectively. $\sym n$
is the symmetric group on $n$ letters, and $\alt n$ is the alternate
subgroup of $\sym n$.  Polynomials and power series are all written
over the formal variable $\hbar$, as is customary in the theory of
quantum algebras.  The Galois field with $p$ elements is written $\F$.
The cyclic group of order $n$ is written $C_n$.

The lower central series of $G$ is $\{\gamma_n(G)\}$, the lower
$p$-central series is $\{P_n(G)\}$, the dimension series is $\{G_n\}$,
the Lie dimension series is $\{L_n(G)\}$, and the derived series is
$G^{(n)}$, and in particular $G'=[G,G]$ --- the definitions shall be
given below.

As is common practice, $H<G$ means that $H$ is a
not-necessarily-proper subgroup of $G$.  For $H<G$, the subgroup of
$H$ generated by $n$-th powers of elements in $H$ is written
$\mho_n(H)$, and $H^{\times n}$ denotes the direct product of $n$
copies of $H$, avoiding the ambiguous ``$H^n$''. The normal closure of
$H$ in $G$ is $H^G$.

Finally, ``$*$'' stands for ``anything'' --- something a speaker would
abbreviate as ``blah, blah, blah'' in a talk. It is used to mean
either that the value is irrelevant to the rest of the computation, or
that it is the only unknown in an equation and therefore does not
warrant a special name.

\subsection{$N$-series}
We first recall a classical construction of Magnus~\cite{magnus:lie},
described for instance in~\cite{lazard:nilp}
and~\cite{huppert-b:fg2}*{Chapter~VIII}.
\begin{defn}
  Let $G$ be a group. An \emdef{$N$-series} is series $\{H_n\}$ of
  normal subgroups with $H_1=G$, $H_{n+1}\le H_n$ and $[H_m,H_n]\leq
  H_{m+n}$ for all $m,n\ge1$.
  The associated Lie ring is
  \[\Lie(G) = \bigoplus_{n=1}^\infty \Lie_n,\]
  with $\Lie_n=H_n/H_{n+1}$ and the bracket operation
  $\Lie_n\otimes\Lie_m\to\Lie_{m+n}$ induced by commutation in $G$.

  For $p$ a prime, an \emdef{$N_p$-series} is an $N$-series
  $\{H_n\}$ such that $\mho_p(H_n)\le H_{pn}$, and the associated Lie ring
  is a restricted Lie algebra over $\F$.
  \[\Lie_{\F}(G) = \bigoplus_{n=1}^\infty \Lie_n,\]
  with the $p$-mapping $\Lie_n\to\Lie_{pn}$ induced by raising to the
  power $p$ in $H_n$.
\end{defn}

We recall that $\Lie$ is a \emph{restricted} Lie algebra
(see~\cite{jacobson:restr} or~\cite{strade-f:mla}*{Section~2.1}) if it
is over a field $\Bbbk$ of characterstic $p$, and there exists a
mapping $x\mapsto x^{[p]}$ such that $\ad x^{[p]}=\ad(x)^p$, $(\alpha
x)^{[p]}=\alpha^px^{[p]}$ and
$(x+y)^{[p]}=x^{[p]}+y^{[p]}+\sum_{i=1}^{p-1}s_i(x,y)$, where the
$s_i$ are obtained by expanding
$\ad(x\otimes\hbar+y\otimes1)^{p-1}(a\otimes1)=\sum_{i=1}^{p-1}s_i(x,y)\otimes
i\hbar^{i-1}$ in $\Lie\otimes\Bbbk[\hbar]$.  Equivalently,
\begin{prop}[Jacobson]\label{prop:jacobson}
  Let $(e_i)$ be a basis of $\Lie$ such that, for some $y_i\in\Lie$,
  we have $\ad(e_i)^p=\ad(y_i)$. Then $\Lie$ is restricted; more
  precisely, there exists a unique $p$-mapping such that
  $e_i^{[p]}=y_i$.
\end{prop}

The standard example of $N$-series is the \emph{lower central series},
$\{\gamma_n(G)\}_{n=1}^\infty$, given by $\gamma_1(G)=G$ and
$\gamma_n(G)=[G,\gamma_{n-1}(G)]$, or the \emph{lower exponent-$p$
  central series} or \emph{Frattini series} given by $P_1(G)=G$ and
$P_n(G)=[G,P_{n-1}(G)]\mho_p(P_{n-1}(G))$. It differs from the lower
central series in that its successive quotients are all elementary
$p$-groups.

The standard example of $N_p$-series is the \emph{dimension series},
also known as the $p$-lower central,
Zassenhaus~\cite{zassenhaus:ordnen}, Jennings~\cite{jennings:gpring},
Lazard~\cite{lazard:nilp} or Brauer series, given by $G_1=G$ and
$G_n=[G,G_{n-1}]\mho_p(G_{\lceil n/p \rceil})$, where $\lceil n/p
\rceil$ is the least integer greater than or equal to $n/p$. It can
alternatively be described, by a result of Lazard~\cite{lazard:nilp},
as
\begin{equation}\label{eq:lazard}
  G_n = \prod_{i\cdot p^j\ge n}\mho_{p^j}(\gamma_i(G)),
\end{equation}
or as
\[G_n = \{g\in G|\,g-1\in\varpi^n\},\]
where $\varpi$ is the augmentation (or fundamental) ideal of the group
algebra $\F G$. Note that this last definition extends to
characteristic $0$, giving a graded Lie algebra $\Lie_\Q(G)$ over
$\Q$. In that case, the subgroup $G_n$ is the isolator of
$\gamma_n(G)$:
\[G_n = \sqrt{\gamma_n(G)} = \{g\in G|\,\langle
g\rangle\cap\gamma_n(G)\neq\{1\}\}.\]
A good reference for these results is~\cite{passi:gr}*{Chapter~VIII}.

We mention finally for completeness another $N_p$-series, the
\emdef{Lie dimension series} $L_n(G)$ defined by
\[L_n(G) = \{g\in G|\,g-1\in\varpi^{(n)}\},\]
where $\varpi^{(n)}$ is the $n$-th Lie power of $\varpi<\Bbbk G$,
given by $\varpi^{(1)}=\varpi$ and
$\varpi^{(n+1)}=[\varpi^{(n)},\varpi]=\{xy-yx|\,x\in\varpi^{(n)},y\in\varpi\}$.
It is then known~\cite{passi-s:liedim} that
\[L_n(G) = \prod_{(i-1)\cdot p^j\ge n}\mho_{p^j}(\gamma_i(G))\]
if $\Bbbk$ is of characteristic $p$, and
\[L_n(G) = \sqrt{\gamma_n(G)}\cap [G,G]\]
if $\Bbbk$ is of characteristic $0$.

In the sequel we will only consider the $N$-series $\{\gamma_n(G)\}$
and $\{P_n(G)\}$ and the $N_p$-series $\{G_n\}$ of dimension
subgroups. We reserve the symbols $\Lie$ and $\Lie_{\F}$ for their
respective Lie algebras.

\begin{defn}
  Let $\{H_n\}$ be an $N$-series for $G$. The \emph{degree} of $g\in
  G$ is the maximal $n\in\N\cup\{\infty\}$ such that $g$ belongs to
  $H_n$.
\end{defn}

A series $\{H_n\}$ has \emdef{finite width} if there is a constant $W$
such that $\ell_n:=\rank[H_n:H_{n+1}]\le W$ holds for all $n$ (Here
$\rank A$ is the minimal number of generators of the abelian group
$A$).  A group has \emdef{finite width} if its lower central series
has finite width --- this definition comes from~\cite{klass-lg-p:fw}.

\begin{defn}
  Let $a=\{a_n\}$ and $b=\{b_n\}$ be two sequences of real numbers. We
  write $a\precsim b$ if there is an integer $C>0$ such that
  $a_n<Cb_{Cn+C}+C$ for all $n\in\N$, and write $a\sim b$ if
  $a\precsim b$ and $b\precsim a$.
\end{defn}
In the sense of this definition, a group has finite width if and only
if $\{\ell_n\}\sim\{1\}$.

I do not know the answer to the following natural
\begin{question}
  If $\rank(\gamma_n(G)/\gamma_{n+1}(G))$ is bounded, does that imply
  that $\rank(G_n/G_{n+1})$, $\rank(P_n(G)/P_{n+1}(G))$ or
  $\rank(L_n(G)/L_{n+1}(G))$ is bounded?  and conversely?
  
  More generally, say an $N$-series $\{H_n\}$ has \emdef{finite width}
  if $\rank(H_n/H_{n+1})$ is bounded over $n\in\N$. If $G$ has a
  finite-width $N$-series intersecting to $\{1\}$, are all $N$-series
  of $G$ of finite width?
\end{question}

The following result is well-known, and shows that sometimes the
Lie ring $\Lie(G)$ is actually a Lie algebra over $\F$.
\begin{lem}\label{lem:liealg}
  Let $G$ be a group generated by a set $S$. Let $\Lie(G)$ be the Lie
  ring associated to the lower central series.
  \begin{enumerate}
  \item If $S$ is finite, then $\Lie_n$ is a finite-rank $\Z$-module
    for all $n$.
  \item If there is a prime $p$ such that all generators $s\in S$ have
    order $p$, then $\Lie_n$ is a vector space over $\F$ for all $n$.
    It then follows that the Frattini series (for that prime $p$) and
    the lower central series coincide.
  \end{enumerate}
\end{lem}
\begin{proof}
  First, $\Lie_1$ is generated by $\overline S$, the image of $S$ in
  $G/G'$. Since $\Lie$ is generated by $\Lie_1$, in particular
  $\Lie_n$ is generated by the finitely many $(n-1)$-fold products of
  elements of $\overline S$; this proves the first point.
  
  Actually, far fewer generators are required for $\Lie_n$; in the
  extremal case when $G$ is a free group, a basis of $\Lie_n$ is given
  in terms of ``standard monomials'' of degree $n$; see
  Subsection~\ref{subs:free} or~\cite{hall:liebasis}.
  
  For the second claim, assume more generally that $s^p\in
  G'$ for all $s\in S$, so that  $G/G'$ is an $\F$-vector space. We
  use the identity $[x,y]^p\equiv [x,y^p]\mod\gamma_3\langle
  x,y\rangle$, due to Philip Hall. Let $g=[x,y]$ be a generator of
  $\gamma_n(G)$, with $x\in G$ and $y\in\gamma_{n-1}(G)$. Then
  $y^p\in\gamma_n(G)$ by induction, so $g^p\in\gamma_{n+1}(G)$ and
  $\Lie_n$ is an $\F$-vector space.
\end{proof}

Anticipating, we note that the groups $\GS$ and $\FG$ we shall
consider satisfy these hypotheses for $p=3$, and $\Gg$ satisfies them
for $p=2$.

\subsection{Growth of groups and vector spaces}
Let $G$ be a group generated by a finite set $S$. The \emdef{length}
$|g|$ of an element $g\in G$ is the minimal number $n$ such that $g$
can be written as $s_1\dots s_n$ with $s_i\in S$. The \emdef{growth
  series} of $G$ is the formal power series
\[\gr(G)=\sum_{g\in G}\hbar^{|g|}=\sum_{n\ge0}f_n\hbar^n,\]
where $f_n=\#\{g\in G|\,|g|=n\}$.  The \emdef{growth function} of $G$ is
the $\sim$-equivalence class of the sequence $\{f_n\}$. Note that
although $\gr(G)$ depends on $S$, this equivalence class is
independent of the choice of $S$.

Let $X$ be a transitive $G$-set and $x_0\in X$ be a fixed base point.
The \emdef{length} $|x|$ of an element $x\in X$ is the minimal length
of a $g\in G$ moving $x_0$ to $x$. The \emdef{growth series} of $X$ is
the formal power series
\[\gr(X,x_0)=\sum_{x\in X}\hbar^{|x|}=\sum_{n\ge0}f_n\hbar^n,\]
where $f_n=\{x\in X|\,\min_{gx_0=x}|g|=n\}$.  The \emdef{growth
  function} of $X$ is the the $\sim$-equivalence class of the sequence
$\{f_n\}$. It is again independent of the choice of $x_0$ and of
generators of $G$.

Let $V=\bigoplus_{n\ge0}V_n$ be a graded vector space. The
\emdef{Hilbert-Poincar\'e} series of $V$ is the formal power series
\[\gr(V)=\sum_{n\ge0}v_n\hbar^n=\sum_{n\ge0}\dim V_n\hbar^n.\]

We return to the dimension series of $G$. Consider the graded algebra
\[\overline{\F G}=\bigoplus_{n=0}^\infty\varpi^n/\varpi^{n+1}.\]
A fundamental result connecting $\Lie_{\F}(G)$ and $\overline{\F G}$ is
the
\begin{thm}[Quillen~\cite{quillen:ab}]\label{thm:quillen}
  $\overline{\F G}$ is the restricted enveloping algebra of the Lie
  algebra $\Lie_{\F}(G)$ associated to the dimension series.
\end{thm}

The Poincar\'e-Birkhoff-Witt Theorem then gives a basis of
$\overline{\F G}$ consisting of monomials over a basis of
$\Lie_{\F}(G)$, with exponents at most $p-1$. As a consequence, we
have the
\begin{prop}[Jennings~\cite{jennings:gpring}]\label{prop:sumprod}
  Let $G$ be a group, and let $\sum_{n\ge1}\ell_n\hbar^n$ be the
  Hilbert-Poincar\'e series of $\Lie_{\F}(G)$. Then
  \[\gr(\overline{\F G})
  =\prod_{n=1}^\infty\left(\frac{1-\hbar^{pn}}{1-\hbar^n}\right)^{\ell_n}.\]
\end{prop}

Approximations from analytical number theory~\cite{li:nt} and complex
analysis give then the
\begin{prop}[\cite{bartholdi-g:lie}, Proposition~2.2
  and~\cite{petrogradsky:polynilpotent}, Theorem~2.1]\label{prop:grlgr}
  Let $G$ be a group and expand the power series
  $\gr(\Lie_{\F}(G))=\sum_{n\ge1}\ell_n\hbar^n$ and
  $\gr(\overline{\F G})=\sum_{n\ge0}f_n\hbar^n$. Then
  \begin{enumerate}
  \item $\{f_n\}$ grows exponentially if and only if $\{\ell_n\}$
    does, and we have
    \[\limsup_{n\to\infty}\frac{\ln \ell_n}n
    =\limsup_{n\to\infty}\frac{\ln f_n}n.\]
  \item If $\ell_n\sim n^d$, then $f_n\sim e^{n^{(d+1)/(d+2)}}$.
  \end{enumerate}
\end{prop}

The Lie algebras we consider have polynomial growth, i.e.\ finite
Gelfand-Kirillov dimension. This notion is more commonly studied for
associative rings~\cite{gelfand-k:dimension}:
\begin{defn}\label{def:gkdim}
  Let $\Lie=\oplus\Lie_n$ be a graded Lie algebra. Its
  \emph{Gelfand-Kirillov dimension} is
  \[\dim_{GK}(\Lie)=\limsup_{n\to\infty}\frac{\log\left(\dim\Lie_1+\dots+\dim{\Lie_n}\right)}{\log n}.\]
\end{defn}
Note that if $\ell_n\sim n^d$, then $\Lie$ has Gelfand-Kirillov
dimension $d+1$. However, the converse is not true, since the sequence
$\log(\ell_1+\dots+\ell_n)/\log n$ need not converge. If the group $G$ has
finite width, then its algebra $\Lie(G)$ has Gelfand-Kirillov
dimension $1$.

Note also that if $A$ is any algebra generated in degree $1$, then
$\dim_{GK}(A)=0$ or $\dim_{GK}(A)\ge1$. Furthermore, George Bergman
showed in~\cite{bergman:growth} that if $A$ is associative, then
$\dim_{GK}(A)=1$ or $\dim_{GK}(A)\ge2$. Victor Petrogradsky showed
in~\cite{petrogradsky:anydim} that there exist Lie algebras of any
Gelfand-Kirillov dimension $\ge1$.

Finally, we recall a connection between the growth of $G$ and that of
$\overline{\F G}$. We use the notation $\sum f_n\hbar^n\ge\sum
g_n\hbar^n$ to mean $f_n\ge g_n$ for all $n\in\N$.

\begin{prop}[\cite{grigorchuk:hp}, Lemma~8]\label{prop:gpgr}
  Let $G$ be a group generated by a finite set $S$. Then
  \[\frac{\gr(G)}{1-\hbar}\ge\gr(\overline{\Bbbk G}).\]
\end{prop}

\section{Branch groups}
Branch groups were introduced by Rostislav Grigorchuk
in~\cite{grigorchuk:jibg}, where he develops a general theory of
groups acting on rooted trees. We shall content ourselves with a
restricted definition; recall that $G\wr\sym d$ is the \emdef{wreath
  product} $G^{\times d}\rtimes\sym d$, the action of $\sym d$ on the
direct product induced by the permutation action of $\sym d$ on
$\Sigma=\{\mathsf1,\dots,\mathsf d\}$.
\begin{defn}
  A group $G$ is \emdef{regular branch} if for some $d\in\N$ there is
  \begin{enumerate}
  \item an embedding $\psi:G\hookrightarrow G\wr\sym d$ such that the
    image of $\psi(G)$ in $\sym d$ acts transitively on $\Sigma$.
    Define for $n\in\N$ the subgroups $\stab_G(n)$ of $G$ by
    $\stab_G(0)=G$, and inductively
    \[\stab_G(n)=\psi^{-1}(\stab_G(n-1)^{\times d})\]
    where $\stab_G(n-1)^{\times d}$ is seen as a subgroup of $G\wr\sym
    d$. One requires then that $\bigcap_{n\in\N}\stab_G(n)=\{1\}$;
  \item a subgroup $K<G$ of finite index with $\psi(K)<K^{\times d}$.
  \end{enumerate}
\end{defn}
To avoid ambiguous bracket notations, we write the decomposition map
\[\psi(g)=\pair<g_1,\dots,g_d>\pi,\]
with $\pi$ expressed as a permutation in disjoint cycle notation.

We shall abbreviate ``regular branch'' to ``branch'', since all the
branch groups that appear in this paper are actually regular branch.
We shall usually omit $d$ from the description, and say that ``$G$
branches over $K$''.

\begin{lem}\label{lem:branchnormal}
  If $G$ is a branch group, then $G$ branches over a subgroup $K$ of
  $G$ such that $K$ is normal in $G$, and $K^{\times d}$ is normal in
  $\psi(K)$.
\end{lem}
\begin{proof}
  Let $G$ be branch over $L$ of finite index, and set $K=\bigcup_{g\in
    G}L^g$, the \emph{core} of $L$. Then obviously $L\triangleleft G$;
  and since $(L^{\times d})^{\psi g}<\psi(K^g)$ for all $g\in G$, we
  have, writing $\psi(g)=\pair<g_1,\dots,g_d>\pi$,
  \[K^{\times d}\le\bigcap_{g\in G}(L^{g_{1^\pi}}\times\dots\times
  L^{g_{d^\pi}})=\bigcap_{g\in G}(L^{\times d})^{\psi g}<K,
  \]
  and $(K^{\times d})^{\psi(g)}=K^{g_{1^\pi}}\times\dots\times
  K^{g_{d^\pi}}=K^{\times d}$, so $K^{\times d}\triangleleft\psi(G)$.
\end{proof}

Let $G$ be a branch group, with $d$, $\Sigma$ and $K$ as in the
definition. The \emdef{rooted tree} on $\Sigma$ is the free monoid
$\Sigma^*$, with root the empty sequence $\emptyset$; it is a metric
space for the distance
\[\dist(\sigma,\tau)=|\sigma|+|\tau|-2\max\{n\in\N|\,\sigma_n=\tau_n\}.\]
The \emdef{natural action} of $G$ is an action on $\Sigma^*$, defined
inductively by
\begin{equation}
  (\sigma_1\sigma_2\dots\sigma_n)^g
  =(\sigma_1)^\pi(\sigma_2\dots\sigma_n)^{g_{\sigma_1}}\text{ for
  }\sigma_1,\dots,\sigma_n\in\Sigma,\label{eq:action}
\end{equation}
where $\psi(g)=\pair<g_1,\dots,g_d>\pi$. By the condition
$\bigcap\stab_G(n)=\{1\}$, this action is faithful and $G$ is
residually finite. Note that $\stab_G(n)$ is the fixator of $\Sigma^n$
in this action.

Note that the action~\eqref{eq:action} gives geometrical meaning to
the branch structure of $G$ that closely parallels the structure of
the tree $\Sigma^*$. Indeed one may consider $G$ as a group acting on
the tree $\Sigma^*$; then the choice of a vertex $\sigma$ of
$\Sigma^*$ and of a subgroup $J$ of $K$ determines a subgroup
$L_\sigma$ of $K$, namely the group of tree-automorphisms of
$\Sigma^*$ that fix $\Sigma^*\setminus\sigma\Sigma^*$ and whose action
on $\sigma\Sigma^*$ is that of an element of $J$ on $\Sigma^*$. The
choice of a subgroup $J_\sigma$ for all $\sigma\in\Sigma^*$ determines
a subgroup $M$ of $K$, namely the closure of the $L_\sigma$ associated
to $\sigma$ and $J_\sigma$ when $\sigma$ ranges over $\Sigma^*$.

This geometrical vision can also give pictorial descriptions of 
the group elements:
\begin{defn}\label{defn:bp}
  Suppose $G$ branches over $K$; let $T$ be a transversal of $K$ in
  $G$, and let $U$ be a transversal of $\psi^{-1}(K^{\times d})$ in
  $K$. The \emdef{branch portrait} of an element $g\in G$ is a
  labeling of $\Sigma^*$, as follows: the root vertex $\emptyset$ is
  labeled by an element of $TU$, and all other vertices are labeled
  by an element of $U$.
  
  Given $g\in G$: write first $g=kt$ with $k\in K$ and $t\in T$; then
  write $k=\psi^{-1}(k_{\mathsf1},\dots,k_{\mathsf d})u_\emptyset$,
  and inductively
  $k_\sigma=\psi^{-1}(k_{\sigma\mathsf1},\dots,k_{\sigma\mathsf
    d})u_\sigma$ for all $\sigma\in\Sigma^*$. Label the root vertex by
  $tu_\emptyset$ and the label the vertex $\sigma\neq\emptyset$ by
  $u_\sigma$.
\end{defn}

There are uncountably many branch portraits, even for a countable
branch group. We therefore introduce the following notion:
\begin{defn}
  Let $G$ be a branch group. Its \emdef{completion} $\overline G$ is
  the inverse limit
  \[\projlim_{n\to\infty}G/\stab_G(n).\]
  This is also the closure in $\aut\Sigma^*$ of $G$ seen through its
  natural action~\eqref{eq:action}.
\end{defn}
Note that since $\overline G$ is closed in $\aut\Sigma^*$ it is a
profinite group, and thus is compact, and totally disconnected. If $G$
has the ``congruence subgroup property''~\cite{grigorchuk:jibg},
meaning that all finite-index subgroups of $G$ contain $\stab_G(n)$
for some $n$, then $\overline G$ is also the profinite completion of
$G$.

\begin{lem}\label{lem:portraits}
  Let $G$ be a branch group and $\overline G$ its completion. Then
  Definition~\ref{defn:bp} yields a bijection between the set of
  branch portraits and $\overline G$.
\end{lem}

We shall often simplify notation by omitting $\psi$ from subgroup
descriptions, as for instance in statements like
``$\stab_G(n)<\stab_G(n-1)^{\times d}$.''

\subsection{The group $\Gg$}
We shall consider more carefully three examples of branch groups in
the sequel. The first example of branch group was considered by
Rostislav Grigorchuk in 1980, and appeared innumerably often in
recent mathematics --- the entire chapter~VIII of~\cite{harpe:ggt} is
devoted to it. It is defined as follows: it is a $4$-generated group
$\Gg$ (with generators $a,b,c,d$), its map $\psi$ is given by
\[\psi:\begin{cases}
  \Gg&\hookrightarrow (\Gg\times\Gg)\rtimes\sym 2\\
  a&\mapsto\makebox[7em][l]{$\pair<1,1>(\mathsf1,\mathsf2),$}
  b\mapsto\pair<a,c>,\\
  c&\mapsto\makebox[7em][l]{$\pair<a,d>,$}d\mapsto\pair<1,b>
\end{cases}\]
and its subgroup $K$ is the normal closure of $[a,b]$, of index $16$.
Rostislav Grigorchuk proved
in~\cites{grigorchuk:burnside,grigorchuk:growth} that $\Gg$ is an
intermediate-growth, infinite torsion group. Its lower central series
was computed in~\cite{bartholdi-g:lie}, along with a description of
its Lie algebra. We shall reproduce that result using a more general
method.

\subsection{The group $\GS$}
This $2$-generated group was introduced by Narain Gupta and Said Sidki
in~\cite{gupta-s:burnside}, where they proved it to be an infinite
torsion group. Later Said Sidki obtained a complete description of its
automorphism group~\cite{sidki:subgroups}, along with information on
its subgroups. It is a branch group with generators $a,t$, its
map $\psi$ is given by
\[\psi:\begin{cases}
  \GS&\hookrightarrow(\GS\times\GS\times\GS)\rtimes\alt3\\
  a&\mapsto\pair<1,1,1>(\mathsf1,\mathsf2,\mathsf3)\\
  t&\mapsto\pair<a,a^{-1},t>,
\end{cases}\]
and its subgroup $K$ is $\GS'$, of index $9$.

The author proved recently~\cite{bartholdi:phd} that $\GS$ has
intermediate growth, which increases its analogy with the Grigorchuk
group mentioned above. An outstanding question was whether $\GS$
has finite width. Ana Cristina Vieira computed
in~\cites{vieira:lcs,vieira:subgroups} the first $9$ terms of the lower
central series and showed that there are all of rank at most $2$. We
shall shortly see, however, that $\GS$ has unbounded width.

The following lemma is straightforward:
\begin{lem}
  $\GS'/(\GS'\times\GS'\times\GS')$ is isomorphic to
 $C_3\times C_3$, generated by $c=[a,t]$ and $u=[a,c]$.
\end{lem}

Note finally that the notations in~\cite{sidki:subgroups} are slightly
different: his $x$ is our $a$, and his $y$ is our $t$.
In~\cite{vieira:lcs} her $y^{[1]}$ is our $u$, and more generally her
$g_1$ is our $\7(g)$ and her $g^{[1]}$ is our $\9(g)$.
In~\cite{bartholdi-g:parabolic}, where a great deal of information on
$\GS$ is gathered, the group is called $\doverline\Gamma$.

\subsection{The group $\FG$}
This other group is at first sight close to $\GS$: it is also
branch, and generated by two elements $a,t$. Its map $\psi$ is given
by
\[\psi:\begin{cases}
  \FG&\hookrightarrow(\FG\times\FG\times\FG)\rtimes\alt3\\
  a&\mapsto\pair<1,1,1>(\mathsf1,\mathsf2,\mathsf3)\\
  t&\mapsto\pair<a,1,t>,
\end{cases}\]
and its subgroup $K$ is $\FG'$, of index $9$.

This group was first considered by Jacek Fabrykowski and Narain
Gupta~\cite{fabrykowski-g:growth2}, who studied its growth.
In~\cite{bartholdi-g:parabolic}, Rostislav Grigorchuk and the author
proved that it is a branch group, and that its subgroup $L=\langle
at,ta\rangle$ has index $3$ and is torsion-free.
In~\cite{bartholdi:phd} another proof of $\FG$'s subexponential growth
is given.

\section{Lie algebras}
We shall now describe the Lie algebras associated to the groups $\Gg$,
$\GS$ and $\FG$ defined in the previous section. We start by
considering a group $G$, and make the following hypotheses on $G$,
which will be satisfied by $\Gg$, $\GS$ and $\FG$:
\begin{enumerate}
\item $G$ is finitely generated by a set $S$;
\item there is a prime $p$ such that all $s\in S$ have order $p$.
\end{enumerate}
Under these conditions, it follows from Lemma~\ref{lem:liealg} that
$\gamma_n(G)/\gamma_{n+1}(G)$ is a finite-dimensional vector space
over $\F$, and therefore that $\Lie(G)$ is a Lie algebra over $\F$
that is finite at each dimension. Clearly the same property holds for
the restricted algebra $\Lie_{\F}(G)$.

We propose the following notation for such algebras:
\begin{defn}\label{defn:liegf}
  Let $\Lie=\bigoplus_{n\ge1}\Lie_n$ be a graded Lie algebra over
  $\F$, and choose a basis $B_n$ of $\Lie_n$ for all $n\ge1$. For
  $x\in\Lie_n$ and $b\in B_n$ denote by $\langle x|b\rangle$ the
  $b$-coefficient of $x$ in base $B_n$.
  
  The \emdef{Lie graph} associated to these choices is an abstract
  graph. Its vertex set is $\bigcup_{n\ge1}B_n$, and each vertex $x\in
  B_n$ has a degree, $n=\deg x$. Its edges are labeled as $\alpha x$,
  with $x\in B_1$ and $\alpha\in\F$, and may only connect a vertex
  of degree $n$ to a vertex of degree $n+1$. For all $x\in B_1$, $y\in
  B_n$ and $z\in B_{n+1}$, there is an edge labeled $\langle
  [x,y]|z\rangle x$ from $y$ to $z$.
  
  If $\Lie$ is a restricted algebra of $\F$, there are additional
  edges, labeled $\alpha\cdot p$ with $\alpha\in\F$, from vertices
  of degree $n$ to vertices of degree $pn$. For all $x\in B_n$ and
  $y\in B_{pn}$, there is an edge labeled $\langle x^p|y\rangle\cdot
  p$ from $x$ to $y$.

  Edges labeled $0x$ are naturally omitted, and edges labeled $1x$
  are simply written $x$.
\end{defn}

There is some analogy between this definition and that of a Cayley
graph --- this topic will be developed in Section~\ref{sec:parabolic}.
The generators (in the Cayley sense) are simply chosen to be the
$\ad(x)$ with $x$ running through $B_1$, a basis of $G/[G,G]$.

A presentation for the $\Lie$ can also be read off its Lie graph. For
every $n$, consider the set $\mathcal W$ of all words of length $n$
over $B_1$. For a path $\pi$ in the Lie graph, define its weight as
the product of the labels on its edges. Each $w\in\mathcal W$ defines
an element of $\Lie_n$, by summing the weights of all paths labeled
$w$ in the Lie graph. Let $\mathcal R_n$ be the set of all linear
dependence relations among these words. Then $\Lie$ admits a
presentation by generators and relations as
\[\Lie = \langle B_1|\,\mathcal R_1,\mathcal R_2,\dots\rangle.\]

Let us give a few examples of Lie graphs. First, if $G$ is abelian,
then its Lie graph has $\rank(G)$ vertices of weight $1$ and no other
vertices. If $G$ is the quaternion group $Q_8=\{\pm1,\pm i,\pm j,\pm
k\}$, then its Lie ring is an algebra over $\F[2]$, and the Lie graph
of $\Lie(Q_8)=\Lie_{\F[2]}(Q_8)$ is
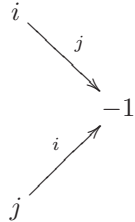
\[\xymatrix{{i}\ar[dr]^{j}\\ & {-1}\\ {j}\ar[ur]^{i}}\]

\subsection{The infinite dihedral group} As another example, let $G$
be the infinite dihedral group $D_\infty=\langle
a,b|\,a^2,b^2\rangle$. Then $\gamma_n(G)=\langle
(ab)^{2^{n-1}}\rangle$ for all $n\ge2$, and its Lie ring is again a
Lie algebra over $\F[2]$, with Lie graph
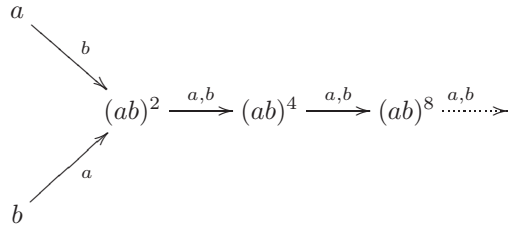
\[\xymatrix{{a}\ar[dr]^{b}\\
  & {(ab)^2}\ar[r]^{a,b} & {(ab)^4}\ar[r]^{a,b} &
  {(ab)^8}\ar@{.>}[r]^{a,b} & {}\\
  {b}\ar[ur]_{a}}\]

Note that the lower $2$-central series of $G$ is different: we have
$G_{2^n}=G_{2^n+1}=\dots=G_{2^{n+1}-1}=\gamma_{n+1}(G)$, so the
Lie graph of $\Lie_{\F[2]}(G)$ is
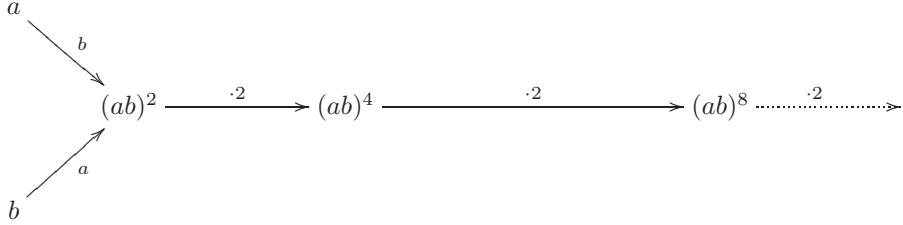
\[\xymatrix{{a}\ar[dr]^{b}\\
  & {(ab)^2}\ar[rr]^{\cdot2} & & {(ab)^4}\ar[rrrr]^{\cdot2} & & & &
  {(ab)^8}\ar@{.>}[rr]^{\cdot2} & & {}\\
  {b}\ar[ur]_{a}}\]

\subsection{The free group}\label{subs:free}
Consider, as an example producing exponential growth, the free group
$F_r$ and its Lie algebra $\Lie$; this is a free Lie algebra of rank
$r$. Using Theorem~\ref{thm:quillen} and M\"obius inversion, we get
\[\dim_\Q(\gamma_n(F_r)/\gamma_{n+1}(F_r)\otimes\Q)=\#\{u\in\mathcal
M|\,\deg u=n\} =\frac1n\sum_{d|n}\mu_{n/d}r^d\precsim r^n,
\]
where $\mu$ is the M\"obius function; therefore $\gr(\overline{\Q
  F_r})\le\frac1{1-r\hbar}$.  Recall that
$\gr(F_r)=\frac{1+\hbar}{1-(2r-1)\hbar}$, so the group growth rate
can be strictly larger than the algebra growth rate in
Proposition~\ref{prop:gpgr}.

It is an altogether different story to find explicitly a basis of
$\Lie$. Pick a basis $X$ of $F_r$; its image in $\Lie_1\cong\Z^r$ is a
generating set of $\Lie$, still written $X$. A \emph{Hall set} is a
linearly ordered set of non-associative words $\mathcal M$ with
$X\subset\mathcal M$ and
\[[u,v]\in\mathcal M\text{ if and only if }u<v\in\mathcal M\text{ and
}(u\in X\text{ or }u=[p,q],q\ge v);\]
furthermore one requires $[u,v]<v$. Note that an order on the
non-associative words uniquely defines a corresponding Hall set.

There are many examples of Hall sets, and for each Hall set $\mathcal
M$ the set $\{u\in\mathcal M|\,|u|=r\}$ is a basis of the abelian
group $\gamma_n(F_r)/\gamma_{n+1}(F_r)$. For example, the \emph{Hall
  basis}~\cite{hall:liebasis} is the linearly ordered set $\mathcal M$
having as maximal elements $X$ in an arbitrary order, and such that
$u<v$ in $\mathcal M$ whenever $\deg(u)>\deg(v)$. It contains then all
$[x,y]$ with $x,y\in X$ and $x>y$; then all $[[u,v],w]$ whenever
$[u,v]<w\le v$ and $u,v,w\in\mathcal M$.

Another basis, more computationally efficient (it is a Lie algebra
equivalent of ``Gr\"obner bases''), is the ``Lyndon-Shirshov
basis''~\cites{shirshov:liebases,lothaire:mots,reutenauer:fla}. It is
defined as follows: order $X$ arbitrarily; on the free monoid $X^*$
put the lexicographical ordering: $u\le uv$, and $uxv<uyw$ for all
$u,v,w\in X^*$ and $x<y\in X$. A non-empty word $w\in X^*$ is a
\emph{Lyndon-Shirshov word} if for any non-trivial factorization
$w=uv$ we have $w<v$. If furthermore we insist that $v$ be
$<$-minimal, then $u$ and $v$ are again Lyndon-Shirshov words. For a
Lyndon-Shirshov word $w$, define its \emph{bracketing} $B(w)$
inductively as follows: if $w\in X$ then $B(w)=w$. If $w=uv$ with $v$
minimal then $B(w)=[B(u),B(v)]$.  Then $\{B(w)\}$ is a basis of
$\Lie$.

From our perspective, an optimal basis $B$ would consist only of
left-ordered commutators, and be prefix-closed, i.e.\ be such that
$[u,x]\in B$ implies $u\in B$; then indeed the Lie algebra structure
of an arbitrary Lie algebra would be determined $\ad(u)$ for all $u\in
B$, and therefore would be a tree in the case of a free Lie algebra.
Kukin announced in~\cite{kukin:liebases} a construction of such bases,
but his proof does not appear to be altogether
complete~\cite{blessenohl-l:liebases}, and the problem of construction
of a left-ordered basis seems to be considered open.

\subsection{The lamplighter group} As another example, consider the
``lamplighter group'' $G=C_2\wr\Z$, with $a$ generating $C_2$ and
$t$ generating $\Z$.  Define the elements
\[a_n=\prod_{i=0}^{n-1}a^{(-1)^i\binom{n-1}it^i}
=at^{-1}a^{-(n-1)}t^{-1}\dots a^{(-1)^{n-1}}t^{n-1}\]
of $G$. Then its Lie algebra $\Lie_{\F[2]}(G)$ is as follows:
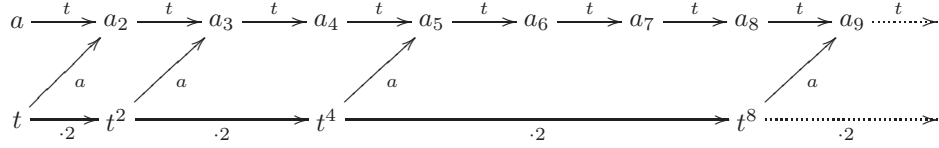
\[\xymatrix{{a}\ar[r]^{t} & {a_2}\ar[r]^{t} & {a_3}\ar[r]^{t} &
  {a_4}\ar[r]^{t} & {a_5}\ar[r]^{t} & {a_6}\ar[r]^{t} &
  {a_7}\ar[r]^{t} & {a_8}\ar[r]^{t} & {a_9}\ar@{.>}[r]^{t} & {}\\
  {t}\ar[ur]_{a}\ar[r]_{\cdot2} &
  {t^2}\ar[ur]_{a}\ar[rr]_{\cdot2} & &
  {t^4}\ar[ur]_{a}\ar[rrrr]_{\cdot2} & & & &
  {t^8}\ar[ur]_{a}\ar@{.>}[rr]_{\cdot2} & & {}}\]
Note that $\Lie_{\F[2]}(G)$ has bounded width, while $G$ has
exponential growth! This shows that in Proposition~\ref{prop:gpgr} the
group growth rate can be exponential while the algebra growth rate is
polynomial.

\subsection{The Nottingham group} As a final example, we give the Lie
graph of the Nottingham group's Lie
algebra~\cites{jennings:sgfps,camina:ng}. Recall that for odd prime $p$
the Nottingham group $J(p)$ is the group of all formal power series
\[\hbar+\sum_{i>1}a_i\hbar^i\in{\F}[[\hbar]],\]
with composition (i.e.\ substitution) as binary operation.
The lower central series is given by
\[J_n=\{\hbar+\sum_{i>\left\lceil\frac{np-1}{p-1}\right\rceil}a_i\hbar^i\},\]
and a basis of $\Lie$ is $\{f_i=\hbar(1+\hbar^i)\}_{i\ge1}$, where
$f_i$ has degree $\left\lfloor\frac{(p-1)i+1}{p}\right\rfloor$. As
basis of $J_1/J_2$, we take
$B_1=\{x=\hbar+\hbar^2+\hbar^3,y=\hbar+\hbar^3\}$. The commutations
are given by
\[[f_i,x]=(i-1)f_{i+1},\qquad
[f_i,y]=\begin{cases}-2f_{i+2}&\text{ if }i\equiv0\mod p\\
  -f_{i+2}&\text{ if }i\equiv1\mod p\\
  0&\text{ otherwise,}
\end{cases}.\]
This gives the Lie graph with ``diamond'' structure~\cite{caranti:nottingham}
{\tiny
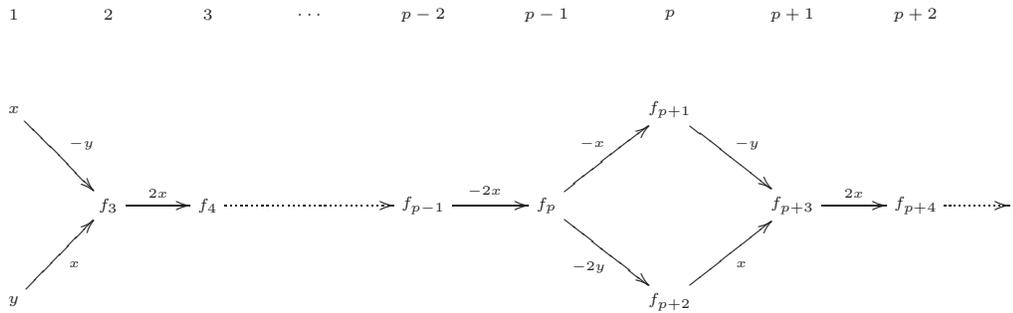
\[\xymatrix{
    {1} & {2} & {3} & {\cdots} & {p-2} & {p-1} & {p} & {p+1} & {p+2}\\
    {x}\ar[dr]^{-y} & & & & & & {f_{p+1}}\ar[dr]^{-y}\\
    & {f_3}\ar[r]^{2x} & {f_4}\ar@{.>}[rr] & & {f_{p-1}}\ar[r]^{-2x} &
    {f_p}\ar[ur]^{-x}\ar[dr]_{-2y} & &
    {f_{p+3}}\ar[r]^{2x} & {f_{p+4}}\ar@{.>}[r] & {}\\
    {y}\ar[ur]_{x} & & & & & & {f_{p+2}}\ar[ur]_{x}}\]}

\subsection{The tree automorphism group's pro-$p$-Sylow
  $\aut_p(\Sigma^*)$}\label{subs:psylow} We start by considering a
typical example of branch group. Let $p$ be prime; write $p'=p-1$ for
notational simplicity. Let $\Sigma$ be the $p$-letter alphabet
$\{\mathsf1,\dots,\mathsf p\}$, and let $x_n$, for $n\in\N$, be the
$p$-cycle permuting the first $p$ branches at level $n+1$ in the tree
$\Sigma^*$. Therefore $x_0$ acts just below the root vertex, and
$x_{n+1}=\pair<x_n,1,\dots,1>$ for all $n$.

For all $n\in\N$ we define $G_n=\aut_p(\Sigma^*)$ as the group
generated by $\{x_0,\dots,x_{n-1}\}$, and $G=\langle
x_0,x_1,\dots\rangle$. Clearly $G=\injlim G_n$, while its closure is
$\overline G=\projlim G_n$.  Note that $G_n$ is a $p$-Sylow of
$\sym{p^n}$, and $\overline G$ is a pro-$p$-Sylow of
$\aut(\Sigma^*)$.

\begin{lem}
  $G=G\wr C_p$; therefore $G$ is a regular branch group over
  itself.
\end{lem}
\begin{proof}
  The subgroup $\langle x_1,x_2,\dots\rangle$ of $G$ is isomorphic to
  $G$ through $x_i\mapsto x_{i-1}$, and its $p$ conjugates under
  powers of $x_0$ commute, since they act on disjoint subtrees.
\end{proof}

Lev Kaloujnine described in~\cite{kaloujnine:struct} the lower central
series of $G_n$, using his notion of \emph{tableau}. Our purpose here
shall be to describe the Lie algebra of $G_n$ (and therefore $G$ and
$\overline G$) using our more geometric approach. Let us just mention
that in Kaloujnine's theory of tableaux his polynomials
$x_1^{e_1}\dots x_n^{e_n}$ correspond to our $\mathbb e_1\dots\mathbb
e_n(x_0)$.

\begin{lem}\label{lem:poisson}
  For $u,v\in G$ and $X,Y\in\{\7,\dots,\usp\}^n$ we have
  \[[X(u),Y(v)]\equiv(X_1+Y_1-\usp)\dots(X_n+Y_n-\usp)
  ([u,v])^{\prod_{i=1}^n(-1)^{p'-Y_i}\binom{X_i}{p'-Y_i}},
  \]
  modulo terms in $[[X(u),Y(v)],G]$.
\end{lem}
\begin{proof}
  The proof follows by induction, and we may suppose $n=1$ without
  loss of generality. Multiplying by terms in $[[X(u),Y(v)],G]$, we
  may assume $Y(v)$ by some element acting only on the last $Y_1$
  subtrees below the root vertex. Then
  \begin{align*}
    [X(u),Y(v)] &\equiv [\pair<u,\dots,u^{(-1)^{X_1}},1,\dots,1>,
    \pair<1,\dots,1,v,\dots,v^{(-1)^{Y_1}}>]\\
    &= \pair<[u,1],\dots,[u^{(-1)^{p'-Y_1}\binom{p}{p'-Y_1}},v],
    \dots,[u^{(-1)^{X_1}},v^{(-1)^{X_1}\binom{p}{X_1}}],\dots,[1,v]>\\
    &\equiv (X+Y-\usp)([u,v])^{(-1)^{p'-Y_1}\binom{p}{p'-Y_1}}.
  \end{align*}
\end{proof}

Note that in the Kaloujnine terminology there is a beautiful
description of $[X(u),Y(v)]$ in terms of Poisson brackets, due to
Vitaly Sushchansky, and due to appear in a forthcoming paper of his.

\begin{thm}\label{thm:psylow}
  Consider the following Lie graph: its vertices are the symbols $X$
  for all words $X\in\{\7,\dots,\usp\}^*$, including the empty
  word $\lambda$.  Their degrees are given by
  \[\deg X_1\dots X_n = 1 + \sum_{i=1}^n X_ip^{i-1}.\]
  For all $m>n\ge0$ and all choices of $X_i$, there is an arrow
  labeled $\7^n$ from $\usp^nX_{n+1}\dots X_m$ to
  $\7^n(X_{n+1}+\8)X_{n+2}\dots X_m$, and an arrow labeled $\7^m$ from
  $\usp^n$ to $\7^n\8\7^{m-n-1}$.
  
  Then the resulting graph is the Lie graph of $\Lie(G)$ and of
  $\Lie_{\F}(G)$.
  
  The subgraph spanned by all words of length up to $n-1$ is the Lie
  graph of $\Lie(G_n)$ and of $\Lie_{\F}(G_n)$.
\end{thm}
\begin{proof}
  We interpret $X$ in the Lie graph as $X(x_0)$ in $G$. The generator
  $x_n$ is then $\7^n(x_0)$. By Lemma~\ref{lem:poisson}, the adjoint
  operators $\ad(x_n)$ correspond to the arrows labeled $\7^n$. The
  arrows connect elements whose degree differ by $1$, so the degree of
  the element $X(x_0)$ is $\deg(X)$ as claimed.

  The power maps $g\mapsto g^p$ are all trivial on the elements
  $X(x_0)$, so the Lie algebra and restricted Lie algebra coincide.

  The elements $X(x_0)$ for $|X|\ge n$ belong to $\stab_G(n)$, and
  hence are trivial in $G_n$.
\end{proof}

\begin{figure}\label{fig:lie:psylow}
  \tiny\[\xymatrix@-2ex{
    {}\ar@{--}[d] &
    {}\ar@{--}[d] &
    {}\ar@{--}[d] &
    {}\ar@{--}[d] &
    {}\ar@{--}[d] &
    {}\ar@{--}[d] &
    {}\ar@{--}[d] &
    {}\ar@{--}[d] \\
    {\7^n}\ar@{--}[d]\ar[r]^{\lambda} &
    {\8\7^{n-1}}\ar@{--}[d]\ar@{.>}[r]^{\lambda} &
    {\usp\7^{n-1}}\ar@{--}[d]\ar[r]^{\7} &
    {\7\8\7^{n-2}}\ar@{--}[d]\ar[r]^{\lambda} &
    {\8\8\7^{n-2}}\ar@{--}[d]\ar@{.>}[r]^{\lambda} &
    {\usp\8\7^{n-2}}\ar@{--}[d]\ar@{.>}[r]^{\lambda} &
    {\usp^2\7^{n-2}}\ar@{--}[d]\ar[r]^{\7^2} &
    {\7^2\8\7^{n-3}}\ar@{--}[d]\ar@{.>}[r] & {} \\
    {\7^3}\ar[r]^{\lambda} & {\8\7\7}\ar@{.>}[r]^{\lambda} &
    {\usp\7\7}\ar[r]^{\7} & {\7\8\7}\ar[r]^{\lambda} &
    {\8\8\7}\ar@{.>}[r]^{\lambda} & {\usp\8\7}\ar@{.>}[r]^{\lambda} &
    {\usp\usp\7}\ar[r]^(0.4){\7^2} & {\7\7\8}\ar@{.>}[r] & {} \\
    {\7^2}\ar[r]^(0.3){\lambda} & {\8\7}\ar@{.>}[r]^{\lambda} &
    {\usp\7}\ar[r]^(0.25){\7} & {\7\8}\ar[r]^{\lambda} &
    {\8\8}\ar@{.>}[r]^{\lambda} & {\usp\8}\ar@{.>}[r]^{\lambda} &
    {\usp\usp}\ar@{.>}[uur]_(0.7){\7^n}\ar[ur]_{\7^3}\\
    {\7}\ar[r]^(0.45){\lambda} & {\8}\ar@{.>}[r]^{\lambda} &
    {\usp}\ar@{.>}[uuur]^{\7^n}\ar[uur]_(0.7){\7^3}\ar[ur]_{\7^2}\\
    {\lambda}\ar@{.>}[uuuur]^(0.6){\7^n}\ar[uuur]_(0.8){\7^3}
    \ar[uur]_(0.7){\7^2}\ar[ur]_{\7}}\]
  \caption{The beginning of the Lie graph of $\Lie(G)$ for $G$ the
    $p$-Sylow of $\aut(\Sigma^*)$}
\end{figure}
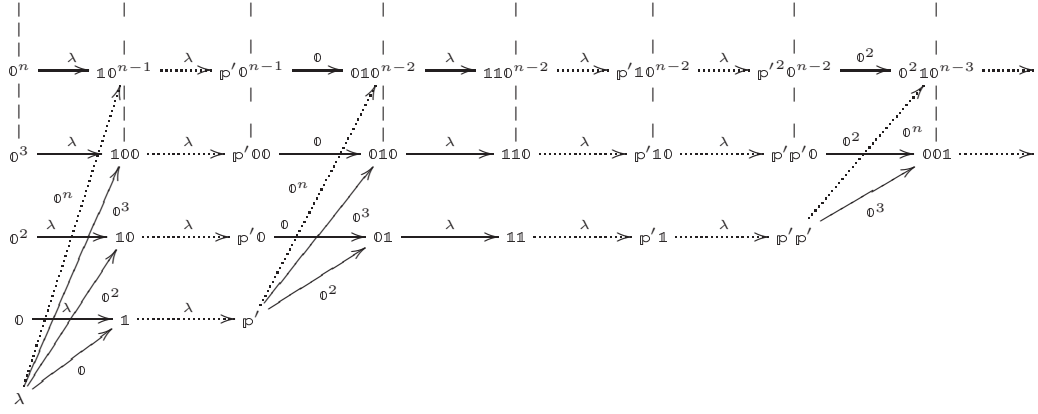

\subsection{The group $\Gg$}\label{subs:Gg}
We give an explicit description of the Lie algebra of $\Gg$, and
compute its Hilbert-Poincar\'e series. These results were obtained
in~\cite{bartholdi-g:lie}, and partly before in~\cite{rozhkov:lcs}.

Set $x=[a,b]$. Then $\Gg$ is branch over $K=\langle x\rangle^\Gg$, and
$K/(K\times K)$ is cyclic of order $4$, generated by $x$.

\newcommand\sm[2]{{\{\begin{smallmatrix}#1\\#2\end{smallmatrix}\}}}

Extend the generating set of $\Gg$ to a formal alphabet
$S=\big\{a,b,c,d,\sm bc,\sm cd,\sm db\big\}$.  Define the
transformation $\sigma$ on words in $S^*$ by
\[\sigma(a)=a\sm bca,\quad\sigma(b)=d,\quad\sigma(c)=b,\quad\sigma(d)=c,\]
extended to subsets by $\sigma\sm xy=\sm{\sigma x}{\sigma y}$. Note
that for any fixed $g\in G$, all elements $h\in\stab_\Gg(1)$ such that
$\psi(h)=\pair<g,*>$ are obtained by picking a letter from each set in
$\sigma(g)$. This motivates the definition of $S$.

\begin{thm}\label{thm:G:struct}
  Consider the following Lie graph: its vertices are the symbols
  $X(x)$ and $X(x^2)$, for words $X\in\{\7,\8\}^*$. Their degrees are
  given by
  \begin{align*}
    \deg X_1\dots X_n(x) = 1 + \sum_{i=1}^n X_i2^{i-1} + 2^n,\\
    \deg X_1\dots X_n(x^2) = 1 + \sum_{i=1}^n X_i2^{i-1} + 2^{n+1}.\\
  \end{align*}
  There are four additional vertices: $a,b,d$ of degree $1$, and
  $[a,d]$ of degree $2$.
  
  Define the arrows as follows: an arrow labeled $\sm xy$ or
  ``$x,y$'' stands for two arrows, labeled $x$ and $y$, and the
  arrows labeled $c$ are there to expose the symmetry of the graph
  (indeed $c=bd$ is not in our chosen basis of $G/[G,G]$):
  \[\xymatrix@R-3ex{
    {a}\ar[r]^-{b,c} & {x} &         {a}\ar[r]^-{c,d} & {[a,d]}\\
    {b}\ar[r]^{a} & {x} &            {d}\ar[r]^-{a} & {[a,d]}\\
    {x}\ar[r]^-{a,b,c} & {x^2} &     {x}\ar[r]^-{c,d} & {\7(x)}\\
    {[a,d]}\ar[r]^-{b,c} & {\7(x)} & {\7*}\ar[r]^-{a} & {\8*}\\
    {\8^n(x)}\ar[r]^{\sigma^n\sm cd} & {\7^{n+1}(x)} &
    {\8^n(x)}\ar[r]^{\sigma^n\sm bd} & {\7^n(x^2)}\\
    {\8^n\7*}\ar[r]^{\sigma^n\sm cd} & {{\7^n\8*}\makebox[0mm][l]{ if
        $n\ge1$.}}}\]
  
  Then the resulting graph is the Lie graph of $\Lie(\Gg)$. A slight
  modification gives the Lie graph of $\Lie_{\F[2]}(\Gg)$: the degree
  of $X_1\dots X_n(x^2)$ is then $2\deg X_1\dots X_n(x)$; and the
  $2$-mappings are given by
  \begin{align*}
    X(x)&\overset{\cdot2}\longrightarrow X(x^2),\\
    \8^n(x^2)&\overset{\cdot2}\longrightarrow \8^{n+1}(x^2).\\
  \end{align*}

  The subgraph spanned by $a,t$, the $X_1\dots X_i(x)$ for $i\le n-2$
  and the $X_1\dots X_i(x^2)$ for $i\le n-4$ is the Lie graph
  associated to the finite quotient $\Gg/\stab_\Gg(n)$.
\end{thm}

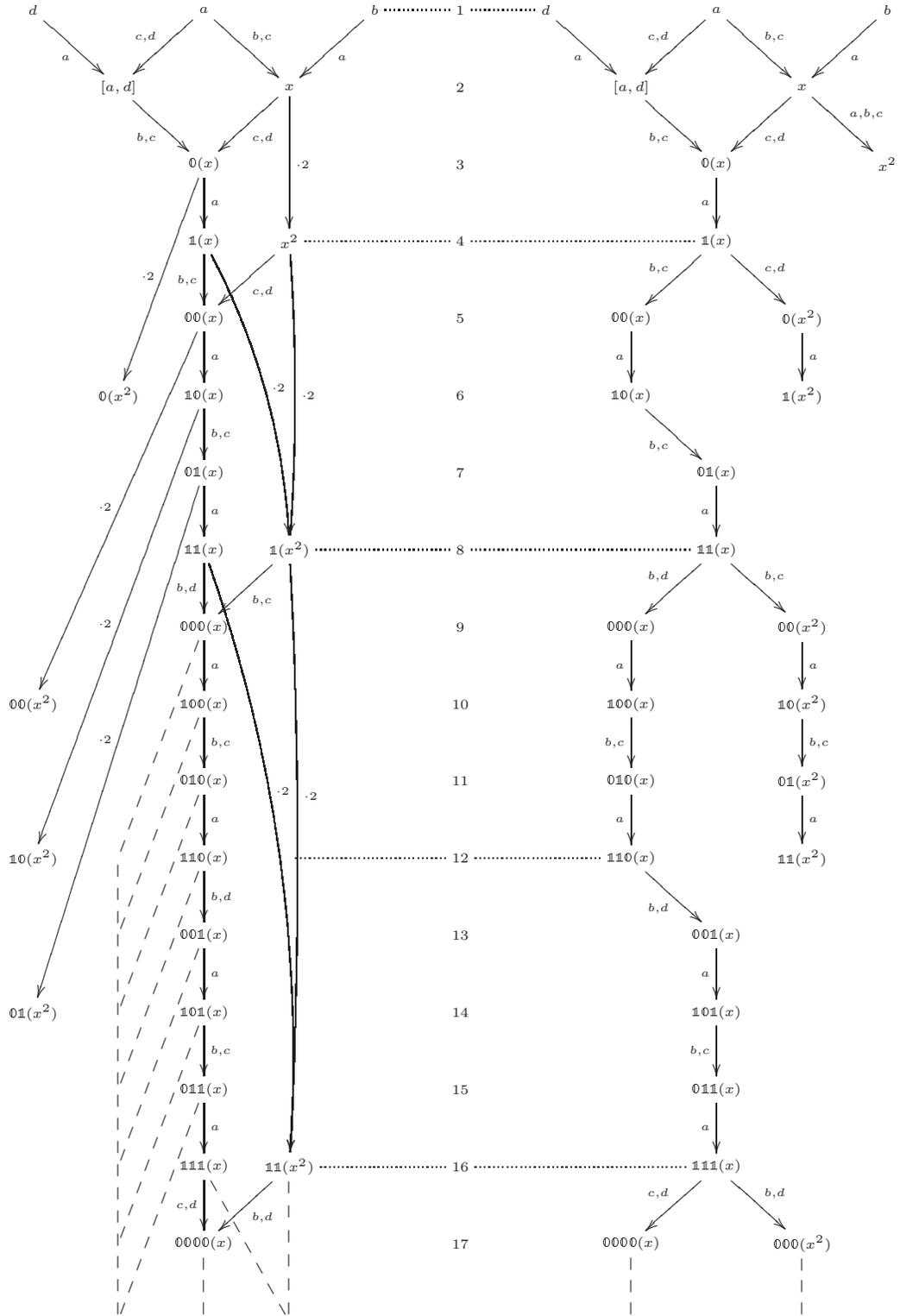
\begin{figure}\label{fig:lie:Gg}
  \tiny\[\xymatrix@d@-7ex@C+2em{
    {b}\ar[dr]^{a} & & {x^2}\\
    & {x}\ar[ur]^{a,b,c}\ar[dr]^{c,d}
    & & & {\7(x^2)}\ar[r]^{a} & {\8(x^2)} & & & {\7\7(x^2)}\ar[r]^{a} &
    {\8\7(x^2)}\ar[r]^{b,c} & {\7\8(x^2)}\ar[r]^{a} & {\8\8(x^2)} & &
    & & & {\7\7\7(x^2)}\ar@{--}[r] &\\
    {a}\ar[ur]^{b,c}\ar[dr]_{c,d}  & & {\7(x)}\ar[r]_{a} &
    {\8(x)}\ar[dr]_{b,c}\ar[ur]^{c,d} & & & {\7\8(x)}\ar[r]_{a} &
    {\8\8(x)}\ar[dr]_{b,d}\ar[ur]^{b,c} & & & & &
    {\7\7\8(x)}\ar[r]_{a} & {\8\7\8(x)}\ar[r]_{b,c} &
    {\7\8\8(x)}\ar[r]_{a} & {\8\8\8(x)}\ar[dr]_{c,d}\ar[ur]^{b,d}\\
    & {[a,d]}\ar[ur]_{b,c} & & & {\7\7(x)}\ar[r]_{a} &
    {\8\7(x)}\ar[ur]_{b,c} & & & {\7\7\7(x)}\ar[r]_{a} &
    {\8\7\7(x)}\ar[r]_{b,c} & {\7\8\7(x)}\ar[r]_{a} &
    {\8\8\7(x)}\ar[ur]_{b,d} & & & & & {\7\7\7\7(x)}\ar@{--}[r] &\\
    {d}\ar[ur]_{a}\\
    {1}\ar@{.}[u]\ar@{.}[d] & {2} & {3} &
    {4}\ar@{.}[uuu]\ar@{.}[dd] &  {5} & {6} & {7} &
    {8}\ar@{.}[uuu]\ar@{.}[dd] & {9} & {10} & {11} &
    {12}\ar@{.}[uu]\ar@{.}[dd] & {13} & {14} & {15} &
    {16}\ar@{.}[uuu]\ar@{.}[dd] & {17}\\
%
    {b}\ar[dr]^{a}\\
    & {x}\ar[dr]^{c,d}\ar[rr]^{\cdot2}
    & & {x^2}\ar[dr]^{c,d}\ar@/^0.2pc/[rrrr]^{\cdot2}
    & & & & {\8(x^2)}\ar[dr]^{b,c}\ar@/^0.3pc/[rrrrrrrr]^(0.4){\cdot2}
    & & & & & & & & {\8\8(x^2)}\ar[dr]^{b,d}\ar@{--}[rr] & &\\
    {a}\ar[ur]^{b,c}\ar[dr]_{c,d} & &
    {\7(x)}\ar[r]^{a}\ar[drrr]_{\cdot2} &
    {\8(x)}\ar[r]_{b,c}\ar@/^0.7pc/[urrrr]^{\cdot2} &
    {\7\7(x)}\ar[r]^{a}\ar[ddrrrrr]_{\cdot2} &
    {\8\7(x)}\ar[r]^{b,c}\ar[ddrrrrrr]_{\cdot2} &
    {\7\8(x)}\ar[r]^{a}\ar[ddrrrrrrr]_{\cdot2} &
    {\8\8(x)}\ar[r]_{b,d}\ar@/^1.2pc/[urrrrrrrr]^(0.4){\cdot2} &
    {\7\7\7(x)}\ar[r]^{a}\ar@{--}[drrr] &
    {\8\7\7(x)}\ar[r]^{b,c}\ar@{--}[drrr] &
    {\7\8\7(x)}\ar[r]^{a}\ar@{--}[drrr] &
    {\8\8\7(x)}\ar[r]^{b,d}\ar@{--}[drrr] &
    {\7\7\8(x)}\ar[r]^{a}\ar@{--}[drrr] &
    {\8\7\8(x)}\ar[r]^{b,c}\ar@{--}[drrr] &
    {\7\8\8(x)}\ar[r]^{a}\ar@{--}[drrr] &
    {\8\8\8(x)}\ar[r]_{c,d}\ar@{--}[urr] &
    {\7\7\7\7(x)}\ar@{--}[r] &\\
    & {[a,d]}\ar[ur]_{b,c} & & & & {\7(x^2)} & & & & & &
    \ar@{--}[rrrrrr] & & & & & & \\
    {d}\ar[ur]_{a} & & & & & & & & &
    {\7\7(x^2)} & & {\8\7(x^2)} & & {\7\8(x^2)}}\]
  \caption{The beginning of the Lie graphs of $\Lie_{\F[2]}(\Gg)$
    (left) and $\Lie(\Gg)$ (right).}
\end{figure}

Figure~\ref{fig:lie:Gg} describes as Lie graphs the top of the Lie
algebras associated to $\Gg$. Note the infinite path, labeled by
\[\sm cda\sigma(\sm cda)\sigma^2(\sm cda)\dots=\sm cda\sm bca\sm
bca\sm bda\sm bca\sm bda\sm bca\sm cda\sm bca\dots;\]
it is the same as the labeling of the parabolic space of $\Gg$ ---
see Section~\ref{sec:parabolic} and~\cite{bartholdi-g:parabolic}.

For proof requires computation, given a term $N$ of a central series
and a generator $s\in\{a,b,c,d\}$, of $[N,x]$ modulo $[N,G,G]$. We
do slightly better in the following lemma --- this will be useful in
Section~\ref{sec:nsGg} where we describe all normal subgroups of $G$.
For that purpose we introduce a symbol $\6(x)=\7(x)\8(x)^{-1}$. We
then have
\[\7(x)=\pair<x,1>,\quad \8(x)=\pair<x,x^{-1}>,\quad \6(x)=\pair<1,x>.\]
\begin{lem}\label{lem:commGg}
  Assume $N$ is a normal subgroup containing the left-hand operand of
  the commutators below. Then modulo $[N,G]'$ we have
  \begin{xalignat*}{2}
    [\7X,a]&=\8X     & [\8X,a]&=\8X^2\\
    [\7X,b]&=\7[X,a] & [\8X,b]&=\7[X,a]+\6[X,c]\\
    [\7X,c]&=\7[X,a] & [\8X,c]&=\7[X,a]+\6[X,d]\\
    [\7X,d]&=1       & [\8X,d]&=\6[X,b]\\
    [x,a]&=x^2       & [x^2,a]&=x^4=\8(x^2+\8x)\\
    [x,b]&=x^2       & [x^2,b]&=\8(x^2+\8x)\\
    [x,c]&=\7(x)+x^2 & [x^2,c]&=\7(x^2+\7x)+\8(x^2+\8x)\\
    [x,d]&=\7(x)     & [x^2,d]&=\7(x^2+\7x)\\
  \end{xalignat*}
\end{lem}
\begin{proof}
  Direct computation, using the decompositions
  $\psi(b)=(a,c)=\7(a)\cdot\6(c)$ etc. and linearizing.
\end{proof}

\begin{proof}[Proof of Theorem~\ref{thm:G:struct}]
  The proof proceeds by induction on length of words, or, what amounts
  to the same, on depth in the lower central series.

  First, the assertion is checked ``manually'' up to degree $3$. The
  details of the computations are the same as
  in~\cite{bartholdi-g:lie}.

  We claim that for all words $X,Y$ with $\deg Y(x)>\deg X(x)$ we have
  $Y(x)\in\langle X(x)\rangle^\Gg$, and similarly $Y(x^2)\in\langle
  X(x^2)\rangle^\Gg$. The claim is verified by induction on $\deg X$.
  
  We then claim that for any non-empty word $X$, either
  $\ad(a)X(*)=0$ (if $X$ starts by ``$\8$'') or
  $\ad(v)X(*)=0$ for $v\in\{b,c,d\}$ (if $X$ starts by ``$\7$'').
  Again this holds by induction.
  
  We then prove that the arrows are as described above; this follow
  from Lemma~\ref{lem:commGg}. For instance,
  \begin{align*}
    \ad(\sigma^n\sm cd)\8^n\7* &= \begin{cases}
      \big(\ad(\sigma^n\sm db)\8^{n-1}\7*,\ad(\sm a1)\8^{n-1}{\7*}\big)\\
      \hfill=\7\ad(\sigma^{n-1}\sm cd)\8^{n-1}\7*=\7^n\8* & \text{ if }n\ge2,\\
      (\ad(\sm bc)\7*,\ad(a)\7*) = \7\8* & \text{ if }n=1.
    \end{cases}
  \end{align*}
  
  Finally we check that the degrees of all basis elements are as
  claimed. For that purpose, we first check that the degree of an
  arrow's destination is always one more than the degree of its
  source.
  Then fix a word $X(*)$, and consider the largest $n$ such that
  $X(*)$ belongs to $\gamma_n(\Gg)$. There is then an expression of
  $X(*)$ as a product of $n$-place commutators on elements of
  $\Gg\setminus[\Gg,\Gg]$, and therefore in the Lie graph there is a
  family of paths starting at some element of $B_1$ and following
  $n-1$ arrows to reach $X(*)$. This implies that the degree of $X(*)$
  is $n$, as required.

  The modification giving the Lie graph of $\Lie_{\F[2]}(\Gg)$ is
  justified by the fact that in $\Lie(\Gg)$ we always have $\deg
  X(x^2)\le2\deg X(x)$, so the element $X(x^2)$ appears always last as
  the image of $X(x)$ through the square map. The degrees are modified
  accordingly. Now $X(x^2)=X\8(x^2)$, and $2\deg X\8(x)\ge 4\deg X(x)$,
  with equality only when $X=\8^n$. This gives an additional square map
  from $\8^n(x^2)$ to $\8^{n+1}(x^2)$, and requires no adjustment of
  the degrees.
\end{proof}

\begin{cor}\label{cor:g:rk}
  Define the polynomials
  \begin{align*}
    Q_2&=-1-\hbar,\\
    Q_3&=\hbar+\hbar^2+\hbar^3,\\
    Q_n(\hbar)&=(1+\hbar)Q_{n-1}(\hbar^2)+\hbar+\hbar^2\text{ for }n\ge4.
  \end{align*}
  Then $Q_n$ is a polynomial of degree $2^{n-1}-1$, and the first
  $2^{n-3}-1$ coefficients of $Q_n$ and $Q_{n+1}$ coincide.  The
  term-wise limit $Q_\infty=\lim_{n\to\infty}Q_n$ therefore exists.
  
  The Hilbert-Poincar\'e series of $\Lie(\Gg/\stab_\Gg(n))$ is
  $3\hbar+\hbar^2+\hbar Q_n$, and the Hilbert-Poincar\'e series of
  $\Lie(\Gg)$ is $3\hbar+\hbar^2+\hbar Q_\infty$.

  The Hilbert-Poincar\'e series of $\Lie_{\F[2]}(\Gg)$ is
  $3+\frac{2\hbar+\hbar^2}{1-\hbar^2}$.

  As a consequence, $\Gg/\stab_\Gg(n)$ is nilpotent of class
  $2^{n-1}$, and $\Gg$ has finite width.
\end{cor}

\begin{proof}
  Consider the sequence of coefficients of $Q_n$. They are, in
  condensed form,
  \[1,2^{2^0},1^{2^0},2^{2^1},1^{2^1},\dots, 2^{n-4},1^{n-4},1^{n-2}.\]
  The $i$th coefficient is $2$ if there are $X(x)$ and $Y(x^2)$ of
  degree $i$ in $\Gg/\stab_\Gg(n)$, and is $1$ if there is only
  $X(x)$. All conclusions follow from this remark.
\end{proof}

\subsection{The group $\GS$}
We now give an explicit description of the Lie algebra of $\GS$,
and compute its Hilbert-Poincar\'e series.

Introduce the following sequence of integers:
\[\alpha_1=1,\quad\alpha_2=2,
\quad\alpha_n=2\alpha_{n-1}+\alpha_{n-2}\text{ for }n\ge3,\]
and $\beta_n=\sum_{i=1}^n\alpha_i$. One has
\begin{align*}
  \alpha_n &= \frac1{2\sqrt2}\left((1+\sqrt2)^n-(1-\sqrt2)^n\right),\\
  \beta_n &= \frac14\left((1+\sqrt2)^{n+1}+(1-\sqrt2)^{n+1}-2\right).
\end{align*}
The first few values are
\[\begin{array}{c|cccccccc}
  n & 1 & 2 & 3 & 4 & 5 & 6 & 7 & 8\\ \hline
  \alpha_n & 1 & 2 & 5 & 12 & 29 & 70 & 169 & 398\\
  \beta_n & 1 & 3 & 8 & 20 & 49 & 119 & 288 & 686
\end{array}\]

\begin{thm}\label{thm:gamma:struct}
  In $\GS$ write $c=[a,t]$ and $u=[a,c]=\9(t)$.  Consider the following
  Lie graph: its vertices are the symbols $X_1\dots X_n(x)$ with
  $X_i\in\{\7,\8,\9\}$ and $x\in\{c,u\}$.  Their degrees are given by
  \begin{align*}
    \deg X_1\dots X_n(c)&=1+\sum_{i=1}^nX_i\alpha_i+\alpha_{n+1},\\
    \deg X_1\dots X_n(u)&=1+\sum_{i=1}^nX_i\alpha_i+2\alpha_{n+1}.
  \end{align*}
  There are two additional vertices, labeled $a$ and $t$, of degree
  $1$.

  Define the arrows as follows:
  \[\xymatrix@R-3ex{
    {a}\tar[r]^{-t} & {c} & {c}\tar[r]^-{t} & {\7(c)}\\
    {t}\aar[r]^{a} & {c} &  {c}\aar[r]^{a} & {u}\\
    & {u}\tar[r]^{t} & {\8(c)}\\
    {\7*}\aar[r]^{a} & {\8*} & {\8*}\aar[r]^{a} & {\9*}\\
    & {\9*}\tar[r]^{t} & {\7\#\makebox[0mm][l]{ whenever
        $*\blue{\overset t\longrightarrow}\#$}}\\
    {\9(c)}\tar[r]^{t} & {\8(u)} &  {\8(c)}\tar[r]^{-t} & {\7(u)}\\
    {\8\7*}\tar[r]^{-t} & {\7\8*} & {\8\8*}\tar[r]^{-t} & {\7\9*}\\
    {\9\7*}\tar[r]^{t} & {\8\8*} &  {\9\8*}\tar[r]^{t} & {\8\9*}}\]
  (Note that these last $3$ lines can be replaced by the rules
  ${\9*}\blue{\overset t\longrightarrow}\8\#$ and
  ${\8*}\blue{\overset{-t}\longrightarrow}\7\#$ for all arrows
  $*\red{\overset a\longrightarrow}\#$.)

  Then the resulting graph is the Lie graph of $\Lie(\GS)$. It is
  also the Lie graph of $\Lie_{\F[3]}(\GS)$, with the only
  non-trivial cube maps given by
  \[\9^n(c)\overset{\cdot3}\longrightarrow\9^n\7\7(c),\qquad
  \9^n(c)\overset{\cdot3}\longrightarrow\9^n\8(u).\]
  
  The subgraph spanned by $a,t$, the $X_1\dots X_i(c)$ for $i\le n-2$
  and the $X_1\dots X_i(u)$ for $i\le n-3$ is the Lie graph associated
  to the finite quotient $\GS/\stab_\GS(n)$.
\end{thm}

\begin{figure}
  \tiny\[\xymatrix@!@-2ex@C-3.4em{
    {1} & {2} & {3} & {4} & {5}\ar@{--}[ddddd] & {6} & {7} & {8} & {9}
    & {10}\ar@{--}[dddd] & {11} & {12} & {13} & {14} & {15}\ar@{--}[d] &
    {16} & {17} & {18} & {19} & {20}\ar@{--}[dd] & {21}\\
    & & & & & & & & & & & & & & {\9\7\7(c)}\tar[dr] & &
    {\7\9\7(c)}\aar[dr] & & {\9\9\7(c)}\tar[dr]\\
    & & & & & & & & & & & & & {\8\7\7(c)}\tar[dr]\aar[ur] & &
    {\8\8\7(c)}\aar[dr]\tar[ur] & & {\8\9\7(c)}\aar[ur] & &
    {\7\8\8(c)}\aar[dr] & & {}\\
    & & & & & & & & & & & & {\7\7\7(c)}\aar[ur] & &
    {\7\8\7(c)}\aar[ur] & & {\9\8\7(c)}\tar[dr]\tar[ur] & &
    {\8\7\8(c)}\tar[ur]\aar[dr] & & {\8\8\8(c)}\aar[ur]\tar[dr]\\
    & & & & & {\7\7(c)}\aar[dr] & & {\9\7(c)}\tar[dr]
    & & {\7\9(c)}\aar[dr] & & {\9\9(c)}\tar[ur]\tar[dr]
    & & & & & & {\7\7\8(c)}\aar[ur] & & {\9\7\8(c)}\tar[ur] & & {}\\
    {a}\tar[dr] & & {\7(c)}\aar[dr] & & {\9(c)}\tar[ur]\tar[dr] & &
    {\8\7(c)}\aar[ur]\tar[dr] & & {\8\8(c)}\tar[ur]\aar[dr] & &
    {\8\9(c)}\aar[ur] & & {\7\8(u)}\aar[dr] & & {\9\8(u)}\tar[dr] & &
    {\9\9(u)}\tar[ur]\\
    & {c}\tar[ur]\aar[dr] & & {\8(c)}\aar[ur]\tar[dr] & &
    {\8(u)}\aar[dr] & & {\7\8(c)}\aar[ur] & &
    {\9\8(c)}\tar[ur]\tar[dr] & & {\8\7(u)}\tar[ur]\aar[dr] & &
    {\8\8(u)}\tar[dr]\aar[ur] & & {\8\9(u)}\aar[ur]\\
    {t}\aar[ur] & & {u}\tar[ur] & & {\7(u)}\aar[ur] & &
    {\9(u)}\tar[ur] & & & & {\7\7(u)}\aar[ur] & & {\9\7(u)}\tar[ur] &
    & {\7\9(u)}\aar[ur]\\
    {2} & {1} & {2} & {1}\ar@{--}[uu] & {2} & {2} & {2} & {2} &
    {1}\ar@{--}[uuu] & {2} & {2} & {2} & {3} & {2} & {4} & {2} & {3} &
    {2} & {2} & {2} & {1}\ar@{--}[uuuuu]}\]
  \caption{The beginning of the Lie graph of $\Lie(\GS)$. The
    generator $\ad(t)$ is shown by plain/blue arrows, and the generator
    $\ad(a)$ is shown by dotted/red arrows.}
\end{figure}

\begin{proof}
  We perform the computations in the completion of $\GS$, still
  written $\GS$. With Lemma~\ref{lem:portraits} in mind, $\GS'$
  is the subgroup generated by all $X(c)$ and $X(u)$, for
  $X\in\{\7,\8,\9\}^*$.
  
  We claim inductively that if $X_i\ge Y_i$ at all positions $i$,
  then $X(c)\in\langle Y(c)\rangle^\GS$, and similarly for $u$.
  Therefore some terms may be neglected in the computations of
  brackets.
  
  Now we compute $\ad(x)y$ for $x,y\in\{a,t,c,u\}$. Here $\equiv$
  means some terms of greater degree have been neglected:
  \begin{align*}
    [a,\7*] &= \8*,\quad [a,\8*] = \9*,
    \quad [a,\9*] = 1\text{ by definition,}\\
    [t,\7*] &= [\pair<a,a^{-1},t>,\pair<*,1,1>]
    =\pair<[a,*],1,1>=\7[a,*]\\
    &\equiv[\pair<a^{-1},t,a>,\pair<*,1,1>]=-\7[a,*],
    \text{ so }[t,\7*]=1,\\
    [t,\8*] &= [\pair<a,a^{-1},t>,\pair<*,*^{-1},1>]\equiv-\7[a,*],\\
    [t,\9*] &= [\pair<a,a^{-1},t>,\pair<*,*,*>]\equiv\8[a,*]+\7[t,*].\\
  \end{align*}
  All asserted arrows follow from these equations.

  Finally, we prove that the degrees of $X(c)$ and $X(u)$ are as
  claimed, by remarking that $\deg c=3$ and $\deg u=4$, that $\deg
  \ad(s)*\ge\deg(*)$ for $s=a,t$ and all words $*$ (so the claimed
  degrees smaller of equal to their actual value), and that each word
  of claimed degree $n$ appears only as $\ad(s)*$ for words $*$ of
  degree at most $n-1$ (so the claimed degrees are greater or equal to
  their actual value).
  
  The last point to check concerns the cube map; we skip the details.
\end{proof}

\begin{cor}\label{cor:gamma:rk}
  Define the following polynomials:
  \begin{align*}
    Q_1&=0,\\
    Q_2&=\hbar+\hbar^2,\\
    Q_3&=\hbar+\hbar^2+2\hbar^3+\hbar^4+\hbar^5,\\
    Q_n&=(1+\hbar^{\alpha_n-\alpha_{n-1}})Q_{n-1}
    +\hbar^{\alpha_{n-1}}(\hbar^{-\alpha_{n-2}}+1+\hbar^{\alpha_{n-2}})Q_{n-2}
    \text{ for }n\ge3.
  \end{align*}
  Then $Q_n$ is a polynomial of degree $\alpha_n$, and the polynomials
  $Q_n$ and $Q_{n+1}$ coincide on their first $2\alpha_{n-1}$ terms.
  The coefficient-wise limit $Q_\infty=\lim_{n\to\infty}Q_n$ therefore
  exists.

  The largest coefficient in $Q_{2n+1}$ is $2^n$, at position
  $\frac12(\alpha_{2n+1}+1)$, so the coefficients of $Q_\infty$ are
  unbounded. The integers $k$ such that $\hbar^k$ has coefficient $1$ in
  $Q_\infty$ are precisely the $\beta_n+1$.
  
  The Hilbert-Poincar\'e series of $\Lie(\GS/\stab_\GS(n))$ is
  $\hbar+Q_n$, and the Hilbert-Poincar\'e series of $\Lie(\GS)$ is
  $\hbar+Q_\infty$. The same holds for the Lie algebra
  $\Lie_{\F[3]}(\GS/\stab_\GS(n))$ and $\Lie_{\F[3]}(\GS)$.

  As a consequence, $\GS/\stab_\GS(n)$ is nilpotent of class
  $\alpha_n$, and $\GS$ does not have finite width.
\end{cor}
\begin{proof}
  Define polynomials
  \[R_n=\sum_{w\in\{\7,\8,\9\}^n}\hbar^{\deg
    w(c)}+\sum_{w\in\{\7,\8,\9\}^{n-1}}\hbar^{\deg w(u)}+\hbar.\]
  Then one checks directly that the polynomials $R_n$ satisfy the same
  initial values and recurrence relation as $Q_n$, hence are
  equal. All convergence properties also follow from the definition of
  $R_n$.
  
  The words of degree $\frac12(\alpha_{2n+1}+1)$ are
  $(\7\8)^{n-1}\7(c)$, $(\7\8)^{n-2}\7\9(u)$, and all the words that
  can be obtained from these by iterating the substitutions
  $\7\7\8\mapsto\8\9\7$, $\8\7\8\mapsto\9\9\7$, $\7\7\9\mapsto\8\9\8$,
  $\8\7\9\mapsto\9\9\8$ along with $\7\8\mapsto\9\7$ and
  $\7\9\mapsto\9\8$ at the beginning of the word. This gives $2^n$
  words in total, half of the form $X(c)$ and half $X(u)$.

  There is a unique word of degree $\beta_n+1$, and that is $\8^n(c)$.
  
  Note that these last two claims have a simple interpretation: there
  are $2^{n-1}$ ways of writing
  $\frac12(\alpha_{2n+1})-1-\alpha_{n+1}$ in base $\alpha$ using only
  the digits $0,1,2$; there is a unique way of writing $\beta_n$ in
  base $\alpha$ using these digits.
\end{proof}

We note as an immediate consequence that
\[[\GS:\gamma_{\beta_n+1}(\GS)]=3^{\frac12(3^n+1)},\]
so that the asymptotic growth of
$\ell_n=\dim(\gamma_n(\GS)/\gamma_{n+1}(\GS))$ is polynomial of degree
$d=\log3/\log(1+\sqrt2)-1$:
\begin{cor}\label{cor:gsgk}
  The Gelfand-Kirillov dimension of $\Lie(\GS)$ is
  $\log3/\log(1+\sqrt2)-1$.
\end{cor}

We then deduce:
\begin{cor}
  The growth of $\GS$ is at least
  $e^{n^{\frac{\log3}{\log(1+\sqrt2)+\log3}}}\cong e^{n^{0.554}}$.
\end{cor}
\begin{proof}
  Apply Proposition~\ref{prop:gpgr} to the series $\sum n^d \hbar^n$,
  which is comparable to the Hilbert-Poincar\'e series of $\Lie(\GS)$.
\end{proof}

Turning to the derived series, we may also improve the general result
$\GS^{(k)}<\gamma_{2^k}(\GS)$ to the following
\begin{thm} For all $k\in\N$ we have
  \[\GS^{(k)}<\gamma_{\alpha_{k+1}}(\GS).\]
\end{thm}
\begin{proof}
  Clearly true for $k=0,1$; then a direct consequence of
  $\GS^{(k)}=\gamma_5(\GS)^{\times3^{k-2}}$ (obtained by Ana
  Vieira in~\cite{vieira:lcs}) and
  $\gamma_{\alpha_j}(\GS)^{\times3}<\gamma_{\alpha_{j+1}}(\GS)$
  for $j=3,\dots,k$.
\end{proof}

\subsection{The group $\FG$}
We now give an explicit description of the Lie algebra of $\FG$,
and compute its Hilbert-Poincar\'e series.

\begin{thm}\label{thm:delta:struct}
  In $\FG$ write $c=[a,t]$ and $u=[a,c]\equiv\9(at)$. For words
  $X=X_1\dots X_n$ with $X_i\in\{\7,\8,\9\}$ define symbols
  $\overline{X_1\dots X_n}(c)$ (representing elements of
  $\FG$) by
  \begin{align*}
    \overline{\mathbb i\7}(c) &= \mathbb i\7(c)/\mathbb i(u),\\
    \overline{\mathbb i\9^{m+1}\8^n}(c) & =\mathbb i\big(
    \overline{\9^{m+1}\8^n}(c)\cdot\7\8^m\7^n(u)^{(-1)^n}\big),\\
    \overline{\mathbb iX}(c) &= \mathbb i\overline X(c)\text{ for all
      other }X.\\
  \end{align*}
  
  Consider the following Lie graph: its vertices are the symbols
  $\overline{X}(c)$ and $X(u)$. Their degrees are given by
  \begin{align*}
    \deg\overline{X_1\dots X_n}(c)&=1+\sum_{i=1}^nX_i3^{i-1}+\frac12(3^n+1),\\
    \deg X_1\dots X_n(u)&=1+\sum_{i=1}^nX_i3^{i-1}+(3^n+1).
  \end{align*}
  There are two additional vertices, labeled $a$ and $t$, of degree
  $1$.

  Define the arrows as follows, for all $n\ge1$:
  \[\xymatrix@R-3ex{
    {a}\tar[r]^{-t} & {c} &      {t}\aar[r]^{a} & {c}\\
    {c}\tar[r]^-{-t} & {\7(c)} & {c}\tar[r]^{a} & {u}\\
    {u}\tar[r]^-{-t} & {\8(c)} & {\overline{\9^n}(c)}\ar[r]^-{-t} &
    {\overline{\7^{n+1}}(c)}\\
    {\7*}\aar[r]^{a} & {\8*} &         {\8*}\aar[r]^{a} & {\9*}\\
    {\9^n\7*}\tar[r]^{t} & {\7^n\8*} & {\9^n\8*}\tar[r]^{t} & {\7^n\9*}\\
    & {\makebox[5mm][r]{$\overline{X_1\dots X_n}(c)$}}
    \tar[r]^{-(-1)^{\sum X_i}t} &
    {\makebox[5mm][l]{$(X_1-1)\dots(X_n-1)(u)$}}}\]
  
  Then the resulting graph is the Lie graph of $\Lie(\FG)$.

  The subgraph spanned by $a,t$, the $\overline{X_1\dots X_i}(c)$ for
  $i\le n-2$ and the $X_1\dots X_i(u)$ for $i\le n-3$ is the Lie graph
  associated to the finite quotient $\FG/\stab_\FG(n)$.
\end{thm}  

\begin{figure}
  \tiny\[\xymatrix@!@-2ex@C-3.4em{
    {1} & {2} & {3} & {4} & {5}\ar@{--}[dd] &
    {6} & {7} & {8} & {9} & {10}\ar@{--}[d] &
    {11} & {12} & {13} & {14} & {15}\ar@{--}[dd] &
    {16} & {17} & {18} & {19} & {20}\ar@{--}[d] & {21}\\
    & & & & & {\overline{\7\7}(c)}\aar[dr] & &
    {\overline{\9\7}(c)}\tar[dr] & &
    {\overline{\8\8}(c)}\aar[dr]\tar[dddr] & &
    {\overline{\7\9}(c)}\aar[dr] & &
    {\overline{\9\9}(c)}\tar[dr]\tar[dddr] & &
    {\overline{\8\7\7}(c)}\aar[dr] & & {\overline{\7\8\7}(c)}\aar[dr]
    & & {\overline{\9\8\7}(c)}\tar[dr]\\
    {a}\tar[dr] & & {\overline{\7}(c)}\aar[dr] & &
    {\overline{\9}(c)}\tar[dr]\tar[ur] & &
    {\overline{\8\7}(c)}\aar[ur] & & {\overline{\7\8}(c)}\aar[ur] & &
    {\overline{\9\8}(c)}\tar[dr]\tar[ur] & &
    {\overline{\8\9}(c)}\aar[ur]\tar[dr] & &
    {\overline{\7\7\7}(c)}\aar[ur] & & {\overline{\9\7\7}(c)}\tar[ur]
    & & {\overline{\8\8\7}(c)}\aar[ur] & & {}\\
    & {c}\tar[ur]\aar[dr] & & {\overline{\8}(c)}\aar[ur]\tar[dr] & &
    {\8(u)}\aar[dr] & & & & & & {\8\7(u)}\aar[dr] & &
    {\7\8(u)}\aar[dr] & & {\9\8(u)}\tar[dr] & & {\8\9(u)}\aar[dr]\\
    {t}\aar[ur] & & {u}\tar[ur] & & {\7(u)}\aar[ur] & & {\9(u)} & & &
    & {\7\7(u)}\aar[ur] & & {\9\7(u)}\tar[ur] & & {\8\8(u)}\aar[ur] &
    & {\7\9(u)}\aar[ur] & & {\9\9(u)}\\
    {2} & {1} & {2} & {1}\ar@{--}[uu] & {2} & {2} & {2} & {1} & {1} &
    {1}\ar@{--}[uuuu] & {2} & {2} & {2} & {2} & {2} & {2} & {2} & {2} &
    {2} & {1} & {1}}\]
  \caption{The beginning of the Lie graph of $\Lie(\FG)$. The
    generator $\ad(t)$ is shown by plain/blue arrows, and the generator
    $\ad(a)$ is shown by dotted/red arrows.}
\end{figure}

\begin{proof}
  The proof is similar to that of Theorems~\ref{thm:G:struct}
  and~\ref{thm:gamma:struct}, but a bit more tricky. Again we perform
  the computations in the completion of $\FG$, still written
  $\FG$. Again $\FG'$ is the subgroup generated by all
  $\overline X(c)$ and $X(u)$, for $X\in\{\7,\8,\9\}^*$.
  
  We claim inductively that if $X_i\ge Y_i$ at all positions $i$,
  then $X(c)\in\langle Y(c)\rangle^\FG$, and similarly for $u$.
  Therefore some terms may be neglected in the computations of
  brackets.
  
  Now we compute $\ad(x)y$ for $x,y\in\{a,t,c,u\}$.  Here $\equiv$
  means some terms of greater degree have been neglected:
  \begin{align*}
    [a,\7*] &= \8*,\quad [a,\8*] = \9*,\quad [a,\9*] = 1\text{ by
      definition,}\\
    [t,\7*] &\equiv [\pair<1,t,a>,\pair<*,1,1>]=1,\\
    [t,\8*] &= [\pair<a,1,t>,\pair<*,*^{-1},1>]=\7[a,*]\\
    &\equiv [\pair<1,t,a>,\pair<*,*^{-1},1>]\equiv -\7[t,*]\\
    [t,\9*] &= [\pair<a,1,t>,\pair<*,*,*>]\equiv \7[a,*]+\7[t,*]+\8[t,*].\\
  \end{align*}
  Note that in the last line the ``negligible'' term $\8[t,*]$ has been
  kept; this is necessary since sometimes the $\7[t,*]$ term cancels
  out.

  Now we check each of the asserted arrows against the
  relations described above. First the ``$a$'' arrows are clearly as
  described, and so are the ``$t$'' arrows on $X(u)$; for instance,
  \begin{align*}
    \ad(t){\9^n\8*}(u) &= \7\ad(a){\9^{n-1}\8*}(u) + \7\ad(t){2^{n-1}\8*}(u)
    + \8\ad(t){\9^{n-1}\8*}(u)\\
    &\equiv \7^n\big(\ad(a){\8*}(u) + \ad(t){\8*}(u)\big) \equiv {\7^n\9*}(u),
  \end{align*}
  which holds by induction on the length of $*$. Next, the ``$t$''
  arrows on $\overline X(c)$ agree; for instance,
  \begin{align*}
    \ad(t)\overline{\9\8^n}(c) &= \7\ad(a)\8^n(c) + \7\ad(t)\8^n(c)
    + \8\ad(t)\8^n(c)\\
    &= \7\9\8^{n-1}(c) + (-1)^n\cdot\7^{n+1}(u) + (-1)^n\cdot\8\7^n(u)\big)\\
    &= \overline{\7\9\8^{n-1}}(c) + (-1)^n\cdot\8\7^n(u)
    \text{ by induction on }n,\\\
    \ad(t)\overline{\9^n}(c) &=
    \ad(t)\9\big(\overline{\9^{n-1}}(c)\cdot\7\8^{n-2}(u)\big)\\
    &\equiv \7\8^{n-1}(u) + \7\big(-\overline{\7^n}(c)-\8^{n-1}(u)\big) +
    \8\big(-\overline{\7^n}(c)-\8^{n-1}(u)\big)\\
    &\equiv-\overline{\7^{n+1}}(c) - \8^n(u).
  \end{align*}
  All other cases are similar. Note how the calculation for
  $\overline{\9\8^n}(c)$ explains the definition of $\overline X(c)$:
  both $\7\9\8^{n-1}(c)$ and $\7^{n+1}(u)$ have degree smaller than
  $d=\deg\overline{\9\8^n}(c)$ in $\Lie(\FG)$, but they are linearly
  dependent in $\gamma_{d-1}(\FG)/\gamma_d(\FG)$.

  Finally, we prove that the degrees of $X(c)$ and $X(u)$ are as
  claimed, by remarking that $\deg c=3$ and $\deg u=4$, that $\deg
  \ad(s)*\ge\deg(*)$ for $s=a,t$ and all words $*$ (so the claimed
  degrees smaller of equal to their actual value), and that each word
  of claimed degree $n$ appears only as $\ad(s)*$ for words $*$ of
  degree at most $n-1$ (so the claimed degrees are greater or equal to
  their actual value).
\end{proof}

\begin{cor}\label{cor:delta:rk}
  Define the integers $\alpha_n=\frac12(5\cdot3^{n-2}+1)$, and the
  polynomials
  \begin{align*}
    Q_2&=1,\\
    Q_3&=1+2\hbar+\hbar^2+\hbar^3+\hbar^4+\hbar^5+\hbar^6,\\
    Q_n(\hbar)&=(1+\hbar+\hbar^2)Q_{n-1}(\hbar^3)+\hbar+
    \hbar^{\alpha_n-2}\text{ for }n\ge4.
  \end{align*}
  Then $Q_n$ is a polynomial of degree $\alpha_n-2$, and the first
  $3^{n-2}+1$ coefficients of $Q_n$ and $Q_{n+1}$ coincide.  The
  term-wise limit $Q_\infty=\lim_{n\to\infty}Q_n$ therefore exists.
  
  The Hilbert-Poincar\'e series of $\Lie(\FG/\stab_\FG(n))$ is
  $2\hbar+\hbar^2Q_n$, and the Hilbert-Poincar\'e series of $\Lie(\FG)$ is
  $2\hbar+\hbar^2Q_\infty$.

  As a consequence, $\FG/\stab_\FG(n)$ is nilpotent of class
  $\alpha_n$, and $\FG$ has finite width.
\end{cor}

\begin{proof}
  Consider the sequence of coefficients of $2\hbar+\hbar^2Q_n$. They
  are, in condensed form,
  \[2,1,2^{3^0},1^{3^0},2^{3^1},1^{3^1},\dots,
  2^{3^{n-3}},1^{3^{n-3}},1^{\frac12(3^{n-1}+1)}.\] The $i$th
  coefficient is $2$ if there are $\overline X(c)$ and $Y(u)$ of
  degree $i$ in $\FG/\stab_\FG(n)$, and is $1$ if there is only
  $\overline X(c)$. All conclusions follow from this remark.
\end{proof}

In quite the same way as for $\GS$, we may improve the general
result $\FG^{(k)}<\gamma_{2^k}(\FG)$:
\begin{thm}
  The derived series of $\FG$ satisfies $\FG'=\gamma_2(\FG)$
  and $\FG^{(k)}=\gamma_5(\FG)^{\times3^{k-2}}$ for $k\ge2$. We
  have for all $k\in\N$
  \[\FG^{(k)}<\gamma_{2+3^{k-1}}(\FG).\]
\end{thm}
\begin{proof}
  It is a general fact for a $2$-generated group $\FG$ that
  $\FG''<\gamma_5(\FG)$. Since $[c,\7(c)]\equiv \7(u)^{-1}$ and
  $[c,u]\equiv \9(c)^{-1}$ (modulo $\gamma_6(\FG)$), we have
  $[\gamma_2(\FG),\gamma_3(\FG)]=\gamma_5(\FG)$ and therefore
  $\FG''=\gamma_5(\FG)$.
  
  Next, $\gamma_5(\FG)=\gamma_3(\FG)^{\times3}\cdot\9(c)$, so
  $\FG^{(3)}=[\gamma_3(\FG),c]^{\times3}=\gamma_5^{\times3}$,
  and the claimed formula holds for all $\FG^{(k)}$ by induction.
  Finally
  $\gamma_{2+3^{j-2}}(\FG)^{\times3}<\gamma_{2+3^{j-1}}(\FG)$
  for all $j=3,\dots,k$.
\end{proof}

We omit altogether the proof of the following two results, since it is
completely analogous to that of Theorem~\ref{thm:delta:struct}.
\begin{thm}
  Keep the notations of Theorem~\ref{thm:delta:struct}. Define now
  furthermore symbols $\overline{X_1\dots X_n}(u)$ (representing
  elements of $\FG$) by
  \begin{align*}
    \overline{\9^n}(u)&=\9^n(u)\cdot \9^{n-1}\7(c)\cdot \9^{n-2}\7\8(c)\cdots
    \9\7\8^{n-2}(c),\\
    \text{and }\overline X(u) &= X(u)\text{ for all other }X.\\
  \end{align*}
  
  Consider the following Lie graph: its vertices are the symbols
  $\overline{X}(c)$ and $\overline X(u)$. Their degrees are given by
  \begin{align*}
    \deg\overline{X_1\dots X_n}(c)&=1+\sum_{i=1}^nX_i3^{i-1}+\frac12(3^n+1),\\
    \deg \9^n(u)&=3^{n+1},\\
    \deg X_1\dots X_n(u)&=\max\{1+\sum_{i=1}^nX_i3^{i-1}+(3^n+1),
    \frac12(9-3^n)+3\sum_{i=1}^nX_i3^{i-1}\}.
  \end{align*}
  There are two additional vertices, labeled $a$ and $t$, of degree
  $1$.

  Define the arrows as follows, for all $n\ge1$:
  \[\xymatrix@R-3ex{
    {a}\tar[r]^{-t} & {c} &      {t}\aar[r]^{a} & {c} & \\
    {c}\tar[r]^-{-t} & {\7(c)} & {c}\tar[r]^{a} & {u} & \\
    {u}\tar[r]^-{-t} & {\8(c)} & {\overline{\9^n}(c)}\ar[r]^{-t} &
    {\overline{\7^{n+1}}(c)}\\
    {\7*}\aar[r]^{a} & {\8*} &         {\8*}\aar[r]^{a} & {\9*}\\
    {\9^n\7*}\tar[r]^{t} & {\7^n\8*} & {\9^n\8*}\tar[r]^{t} & {\7^n\9*}\\
    & {\makebox[5mm][r]{$\overline{X_1\dots X_n}(c)$}}
    \tar[r]^{-(-1)^{\sum X_i}t} &
    {\makebox[5mm][l]{$(X_1-1)\dots(X_n-1)(u)$}}\\
    {c}\ar[r]^{\cdot3} & {\overline{\7\7}(c)} &
    {\overline{2^n}(u)}\ar[r]^{\cdot3} & {\overline{2^{n+1}}(u)}\\
    {*\7(c)}\ar[r]^{\cdot3} & {*\9(u)\makebox[0mm][l]{ if
        $3\deg*\7(c)=\deg*\9(u)$}}}\]
  
  Then the resulting graph is the Lie graph of $\Lie_{\F[3]}(\FG)$.
  
  The subgraph spanned by $a,t$, the $\overline{X_1\dots X_i}(c)$ for
  $i\le n-2$ and the $X_1\dots X_i(u)$ for $i\le n-3$ is the Lie graph
  of the Lie algebra $\Lie_{\F[3]}(\FG/\stab_\FG(n))$.

  As a consequence, the dimension series of $\FG/\stab_\FG(n)$
  has length $3^{n-1}$ (the degree of $\overline{2^n}(u)$), and
  $\FG$ has finite width.
\end{thm}

We have then from Proposition~\ref{prop:grlgr}
\begin{cor}
  The growth of $\FG$ is at least $e^{n^{\frac12}}$.
\end{cor}

\section{Parabolic Space}\label{sec:parabolic}
In the natural action of a branch group $G$ on the tree $\Sigma^*$,
consider a ``parabolic subgroup'' $P$, the stabilizer of an infinite
ray in $\Sigma^*$. (The terminology comes from geometry, where a
parabolic subgroup is the stabilizer of a point on the boundary of an
appropriate $G$-space). Such a parabolic subgroup may be defined
directly as follows: let
$\omega=\omega_1\omega_2\dots\in\Sigma^\infty$ be an infinite
sequence. Set $P_0^\omega=G$ and inductively
\[P_n^\omega = \psi^{-1}(G\times\dots\times P_{n-1}^\omega\dots\times G),\]
with the `$P_{n-1}^\omega$' in position $\omega_n$. Set
$P^\omega=\bigcap_{n\ge0}P_n^\omega$.

In the natural tree action~\eqref{eq:action} of $G$ on $\Sigma^*$ or
on $\Sigma^\infty$ its boundary, $P_n^\omega$ is the stabilizer of the
point $\omega_1\dots\omega_n$, and $P^\omega$ is the stabilizer of the
infinite sequence $\omega$.

The following facts easily follow from the definitions:
\begin{lem}
  $\bigcap_{\omega\in\Sigma^\infty}P^\omega = 1$. The index of
  $P_n^\omega$ in $G$ is $d^n$, and that of $P^\omega$ is infinite.
\end{lem}

\begin{defn}
  Let $G$ be a branch group. A \emdef{parabolic space} for $G$ is a
  homogeneous space $G/P$, where $P$ is a parabolic subgroup.
\end{defn}

Suppose now that $G$ is finitely generated by a set $S$.
\begin{prop}[\cite{bartholdi-g:spectrum}]\label{prop:polygrowth}
  Suppose that the length $|\cdot|$ on the branch group $G$ satisfies
  the following condition: there are constants $\lambda,\mu$ such that
  for all $g\in\stab_G(1)$, writing $\psi(g)=(g_1,\dots,g_d)$, one has
  $|g_i|<\lambda|g|+\mu$.

  Then all parabolic spaces of $G$ have polynomial growth of degree at
  most $\log_{1/\lambda}(d)$.
\end{prop}

\begin{thm}\label{thm:growth}
  Let $G$ be a finitely generated branch group. Then there exists a
  constant $C$ such that for any $x_0\in G$ we have
  \[\frac{C\gr(G/P,x_0P)}{1-\hbar}\ge \frac{\gr\Lie(G)}{1-\hbar}.\]
\end{thm}
\begin{proof}
  Assume $G$ acts on a $d$-regular tree, and write as before $d'=d-1$.
  The proof relies on an identification of the Lie action on group
  elements and the natural action on tree levels. We first claim that
  for any $u\in K$ and $W\in\{\6,\dots,\mathbb d'\}^*$
  \[\deg W(u)\ge \deg(\7^{|W|}(u)) + d_{G/P}(\7^{|W|},W),\]
  where $d(W,X)$ is the length of a minimal word moving $W$ to $X$ in
  the tree $\Sigma^*$.
  
  Therefore the growth of $\Lie(G)$ and $G/P$ may be compared just by
  considering the degrees of elements of the form $\7^n(u)$ for some
  fixed $u\in K$; indeed the other $W(u)$ will contribute a smaller
  growth to the Lie growth series than the corresponding vertices to
  the parabolic growth series, and the $N$ finitely many values $u$
  may take in a branch portrait description will be taken care of by
  the constant $C$.

  Now there is a constant $\ell\in\N$ such that $\7^{\ell+m}(u)$ has
  greater degree than $(\mathbb d')^m(u)$ for all $m\in\N$. Indeed
  there exists $k\in K$ and $\ell\in\N$ such that
  $[k,u]=\7^\ell(u)$, and then
  $[\7^mk,{\mathbb d'}^m(u)]=\7^{\ell+m}(u)$, proving the claim.

  We may now take $C=\ell N$. The Lie growth series is the sum over
  all $n\in\N$ and coset representatives $u\in T$ of the
  power series counting the growth of $W(u)$ over words $W$ of length
  $n$. There are $N$ choices for $u$, and for given $u$ at most $\ell$
  of these power series overlap.
%
%
\end{proof}

Note that this result is valid even if the action on the rooted tree
is not cyclic, i.e.\ even if in the decomposition map $G\to G\wr A$
the finite group $A$ is not cyclic. If $A$ is not nilpotent, then the
Lie algebra $\Lie$ is no longer isomorphic to $G$, so the best we can
hope for is an inequality bounding the growth of $\Lie$ by that of $G/P$.

There are examples of groups of exponential growth, whose Lie algebra
has subexponential growth --- see for example~\cite{bartholdi:zp2-1}.

%
%
%

\section{Normal Subgroups}\label{sec:nsGg}
Using the notion of branch portrait, it is not too difficult to
determine the exact structure of normal subgroups in a branch group.
Consider a $p$-group $G$ and its $p$-Lie algebra $\Lie$ over $\F$.
Normal subgroups of $G$ correspond to ideals of $\Lie$, just as
subgroups of $G$ correspond to subalgebras of $\Lie$; and the index of
$H<G$ is $p^{\dim \Lie/\mathcal M}$, where the subgroup $H$
corresponds to the subalgebra $\mathcal M$. This correspondence is not
exact, and we shall neither use it nor make it explicit; however it
serves as a motivation for relating subgroup growth and the study of
Lie algebras. In all cases, sufficient knowledge of $\Lie$, as well as
its finiteness of width, allow an explicit description of the normal
subgroup lattice of $G$.

We focus on the first and most important example, $\Gg$, for which we
obtain an explicit answer. The computations presented here clearly
extend, \emph{mutatis mutandis}, to any regular branch group.

Set $\mathcal W=\{\7,\8\}^*$, and order words $X\in\mathcal W$ by
``reverse shortlex'': the rank of $X_1\dots X_n$ is
\[\#X_1\dots X_n = 1 + \sum_{i=1}^n X_i2^{i-1} + 2^n.\]
(Note that $\#X=\deg X(x)$ according to the definition in
Subsection~\ref{subs:Gg}.) We write $<$ the order induced by rank.

\begin{thm}\label{thm:normal:Gg}
  The non-trivial normal subgroups of $\Gg$ are as follows:
  \begin{itemize}
  \item there are respectively $1,7,7,1$ subgroups of index $1,2,4,8$
    corresponding to the lifts to $\Gg$ of subgroups of
    $\Gg/[\Gg,\Gg]=C_2^{\times3}$;
  \item there are $12$ other subgroups of $\Gg$ not contained in $K$:
    six of index $8$, namely $\langle[a,c],d^ab\rangle^\Gg$, $\langle
    c\rangle^\Gg$, $\langle x,c^ad\rangle^\Gg$, $\langle
    b\rangle^\Gg$, $\langle[a,d],b^ac\rangle^\Gg$, and
    $\langle d,x^2\rangle^\Gg$; four of index $16$, namely
    $\langle[a,c]\rangle^\Gg$, $\langle[a,d],x^2\rangle^\Gg$,
    $\langle d\rangle^\Gg$, and $\langle[a,d],x^2d\rangle^\Gg$;
    and two of index $32$, namely $\langle[a,d]x^2\rangle^\Gg$ and
    $\langle[a,d]\rangle^\Gg$;
  \item all normal subgroups $N\triangleleft\Gg$ contained in $K$ are
    of the form
    \begin{equation}\label{eq:defW}
      W(A;B_1,\dots,B_m;C):=
      \langle A(x)B_1(x^2)\dots B_m(x^2),C(x^2)\rangle^\Gg,
    \end{equation}
    for words $A,B_i,C\in\mathcal W$. Assume the functions
    $M(A,\{B_i\},C)$ and $S(A,\{B_i\},C)$, with values in $\mathcal
    W$, defined in the proof. Then there is a unique description of
    $N$ in the form~\ref{eq:defW} satisfying $B_1<B_2<\dots<B_m\le
    S(A,\{B_i\},C)<C\le M(A,\{B_i\})$.
    
    The index of $N$ is $2^{\#A+\#S(A,\{B_i\},C)}$.

    The groups can furthermore be subdivided in three types:
    \begin{description}
    \item[$\rmI$] $C\le\7^{|A|}$ and $A\le\7^{|C|}\8\7$. Then all $B_i$
      are optional, i.e.\ there are $2^m$ groups with these $A$ and
      $C$, obtained by choosing any subset of the $B_i$'s;
    \item[$\rmII$] $C>\7^{|A|}$ and $C\le\7^{|A|+1}$. Then $A=B_1\8$ and
      all other $B_i$'s are optional;
    \item[$\rmIII$] $A=\7^n$ and some $B_i=\7^{n-1}$. Then in fact an
      alternate description exists, obtained by suppressing $A$ and
      $B_i$ from the description.
    \end{description}
  \end{itemize}
\end{thm}
Note that we have only described finite-index subgroups of $\Gg$.
Since $\Gg$ is just-infinite, all its non-trivial normal subgroups
have finite index.

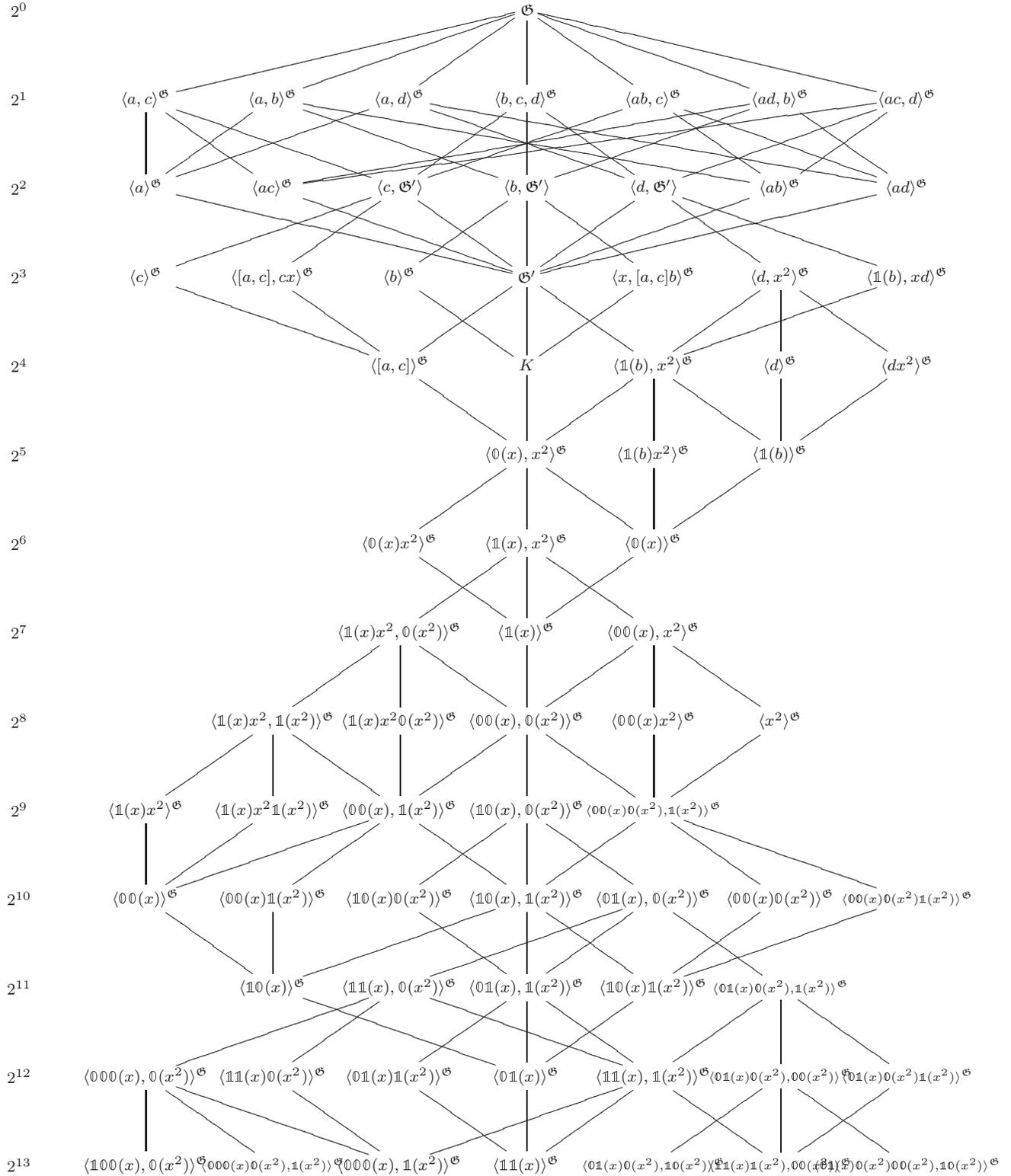
\begin{figure}
  \footnotesize\[\kern-2cm\xymatrix@!@R-3ex@C+1em{
    {2^0} &
    & & & {\Gg}\lar[dlll]\lar[dll]\lar[dl]\lar[d]\lar[dr]\lar[drr]\lar[drrr]\\
    {2^1} &
    \squish{a,c}\lar[d]\lar[dr]\lar[drr] &
    \squish{a,b}\lar[dl]\lar[drr]\lar[drrrr] &
    \squish{a,d}\lar[dll]\lar[drr]\lar[drrrr] &
    \squish{b,c,d}\lar[dl]\lar[d]\lar[dr] &
    \squish{ab,c}\lar[dll]\lar[dr]\lar[drr] &
    \squish{ad,b}\lar[dllll]\lar[dll]\lar[dr] &
    \squish{ac,d}\lar[dlllll]\lar[dll]\lar[dl]\\
    {2^2} &
    \squish{a}\lar[drrr] &
    \squish{ac}\lar[drr] &
    {\langle c,\Gg'\rangle}\lar[dll]\lar[dl]\lar[dr] &
    {\langle b,\Gg'\rangle}\lar[dl]\lar[d]\lar[dr] &
    {\langle d,\Gg'\rangle}\lar[dl]\lar[dr]\lar[drr] &
    \squish{ab}\lar[dll] &
    \squish{ad}\lar[dlll]\\
    {2^3} &
    \squish{c}\lar[drr] &
    \squish{[a,c],cx}\lar[dr] &
    \squish{b}\lar[dr] &
    {\Gg'}\lar[dl]\lar[d]\lar[dr] &
    \squish{x,[a,c]b}\lar[dl] & 
    \squish{d,x^2}\lar[dl]\lar[d]\lar[dr] &
    \squish{\8(b),xd}\lar[dll]\\
    {2^4} &
    & & \squish{[a,c]}\lar[dr]
    & {K}\lar[d]
    & \squish{\8(b),x^2}\lar[dl]\lar[d]\lar[dr]
    & \squish{d}\lar[d]
    & \squish{dx^2}\lar[dl]\\
    {2^5} &
    & & & \squish{\7(x),x^2}\lar[dl]\lar[d]\lar[dr]
    & \squish{\8(b)x^2}\lar[d]
    & \squish{\8(b)}\lar[dl]\\
    {2^6} &
    & & \squish{\7(x)x^2}\lar[dr] &
    \squish{\8(x),x^2}\lar[dl]\lar[d]\lar[dr] &
    \squish{\7(x)}\lar[dl]\\
    {2^7} &
    & & \squish{\8(x)x^2,\7(x^2)}\lar[dl]\lar[d]\lar[dr] &
    \squish{\8(x)}\lar[d] &
    \squish{\7\7(x),x^2}\lar[dl]\lar[d]\lar[dr]\\
    {2^8} &
    & \squish{\8(x)x^2,\8(x^2)}\lar[dl]\lar[d]\lar[dr] &
    \squish{\8(x)x^2\7(x^2)}\lar[d] &
    \squish{\7\7(x),\7(x^2)}\lar[dl]\lar[d]\lar[dr] &
    \squish{\7\7(x)x^2}\lar[d] &
    \squish{x^2}\lar[dl]\\
    {2^9} &
    \squish{\8(x)x^2}\lar[d] &
    \squish{\8(x)x^2\8(x^2)}\lar[dl] &
    \squish{\7\7(x),\8(x^2)}\lar[dll]\lar[dl]\lar[dr] &
    \squish{\8\7(x),\7(x^2)}\lar[dl]\lar[d]\lar[dr] &
    \ssquish{\7\7(x)\7(x^2),\8(x^2)}\lar[dl]\lar[dr]\lar[drr]\\
    {2^{10}} &
    \squish{\7\7(x)}\lar[dr] &
    \squish{\7\7(x)\8(x^2)}\lar[d] &
    \squish{\8\7(x)\7(x^2)}\lar[dr] &
    \squish{\8\7(x),\8(x^2)}\lar[dll]\lar[d]\lar[dr] &
    \squish{\7\8(x),\7(x^2)}\lar[dll]\lar[dl]\lar[dr] &
    \squish{\7\7(x)\7(x^2)}\lar[dl] &
    \ssquish{\7\7(x)\7(x^2)\8(x^2)}\lar[dll]\\
    {2^{11}} &
    & \squish{\8\7(x)}\lar[drr] &
    \squish{\8\8(x),\7(x^2)}\lar[dll]\lar[dl]\lar[drr] &
    \squish{\7\8(x),\8(x^2)}\lar[dl]\lar[d]\lar[dr] &
    \squish{\8\7(x)\8(x^2)}\lar[dl] &
    \ssquish{\7\8(x)\7(x^2),\8(x^2)}\lar[dl]\lar[d]\lar[dr]\\
    {2^{12}} &
    \squish{\7\7\7(x),\7(x^2)}\lar[d]\lar[dr]\lar[drr] &
    \squish{\8\8(x)\7(x^2)}\lar[dr] &
    \squish{\7\8(x)\8(x^2)}\lar[dr] &
    \squish{\7\8(x)}\lar[d] &
    \squish{\8\8(x),\8(x^2)}\lar[dll]\lar[dl]\lar[dr] &
    \ssquish{\7\8(x)\7(x^2),\7\7(x^2)}\lar[dl]\lar[d]\lar[dr] &
    \ssquish{\7\8(x)\7(x^2)\8(x^2)}\lar[dl]\\
    {2^{13}} &
    \squish{\8\7\7(x),\7(x^2)} &
    \ssquish{\7\7\7(x)\7(x^2),\8(x^2)} &
    \squish{\7\7\7(x),\8(x^2)} &
    \squish{\8\8(x)} &
    \ssquish{\7\8(x)\7(x^2),\8\7(x^2)} &
    \ssquish{\8\8(x)\8(x^2),\7\7(x^2)} &
    \ssquish{\7\8(x)\7(x^2)\7\7(x^2),\8\7(x^2)}}\]
  \caption{The top of the lattice of normal subgroups of $\Gg$, of
    index at most $2^{13}$}\label{fig:latticeG}
\end{figure}

We represent the top of the lattice in Figure~\ref{fig:latticeG},
containing all normal subgroups of index at most $2^{13}$ (there are
never more than $7$ subgroups of given lesser index).

The first few subgroups of $K$ are described in Table~\ref{table:ns},
sorted by their index in $\Gg$, and identified by their type in
$\{(\rmI),(\rmII),(\rmIII)\}$. We write $\lambda$ for the empty
sequence. An argument $[B_i]$ means that term is optional, and
therefore stands for two groups, one with that term and one without.

Among the remarkable subgroups are: the subgroup $K^{\times
  2^n}=\langle \7^n(x)\rangle^\Gg$, written $K_n$
in~\cite{bartholdi-g:parabolic}; the subgroup $K^{\times
  2^n}\mho_2(K)^{\times 2^{n-1}}=\langle \7^n(x),\7^{n-1}(x^2)\rangle$,
written $N_n$ in~\cite{bartholdi-g:parabolic}; and $\stab_G(n)=\langle
\7^{n-3}(\8(x)x^2),\7^{n-2}(x^2)\rangle$.

The lattice of normal subgroups of $\Gg$ is described in
Figure~\ref{fig:latticeG}. Even though I do not understand completely
the lattice's structure, some remarks can be made: the lattice has a
fractal appearance; all its nodes have $1$ or $3$ descendants, and $1$
or $3$ ascendants. Large portions of it have a grid-like structure.
This can be explained by the construction $N\rightsquigarrow N\times
N$ of normal subgroups, lending the lattice some self-similarity.

\begin{table}
  \begin{center}
    \begin{tabular}{r|c|lll}
    Index  & Count & Description\\ \hline
    $2^4$&    $1$&      $W(\lambda;;\lambda)_\rmI=K$\\
    $2^5$&    $1$&      $W(\7;;\lambda)_\rmI$\\
    $2^6$&    $3$&      $W(\8;;\lambda)_\rmI$& $W(\7;[\lambda];\7)_\rmI$\\
    $2^7$&    $3$&      $W(\7\7;;\lambda)_\rmI$& $W(\8;[\lambda];\7)_\rmI$\\
    $2^8$&   $5$&      $W(\8\7;;\lambda)_\rmI$& $W(\7\7;;\7)_\rmI$\\
    &&&$W(\8;\lambda,[\7];\8)_\rmII$& $W(\infty;\lambda,\7;\8)_\rmIII$\\
    $2^9$&   $5$&      $W(\8\7;;\7)_\rmI$& $W(\7\7;[\7];\8)_\rmI$&
    $W(\8;\lambda,[\8];\7\7)_\rmII$\\
    $2^{10}$&   $7$&      $W(\7\8;;\7)_\rmI$& $W(\8\7;[\7];\8)_\rmI$&
    $W(\7\7;[\7],[\8];\7\7)_\rmI$\\
    $2^{11}$&  $5$&      $W(\8\8;;\7)_\rmI$& $W(\7\8;[\7];\8)_\rmI$&
    $W(\8\7;[\8];\7\7)_\rmI$\\
    $2^{12}$&  $7$&      $W(\7\7\7;;\7)_\rmI$& $W(\8\8;[\7];\8)_\rmI$&
    $W(\7\8;[\7],[\8];\7\7)_\rmI$\\
    $2^{13}$&  $7$&      $W(\8\7\7;;\7)_\rmI$& $W(\7\7\7;[\7];\8)_\rmI$&
    $W(\8\8;[\8];\7\7)_\rmI$\\
    &&&&$W(\7\8;\7,[\7\7];\8\7)_\rmII$\\
    $2^{14}$&$13$&    $W(\7\8\7;;\7)_\rmI$& $W(\8\7\7;[\7];\8)_\rmI$&
    $W(\7\7\7;;\7\7)_\rmI$\\
    &&&$W(\8\8;\8,[\7\7];\8\7)_\rmII$& $W(\7\8;\7,[\7\7],[\8\7];\7\8)_\rmII$\\
    &&&$W(\infty;\8,0\7;\8\7)_\rmIII$& $W(\infty;\7,[\8],0\7;\7\8)_\rmIII$\\
    $2^{15}$&$9$&     $W(\7\8\7;;\8)_\rmI$& $W(\8\7\7;;\7\7)_\rmI$&
    $W(\7\7\7;[\7\7];\8\7)_\rmI$\\
    &&$W(\8\8;\8,[\8\7];\7\8)_\rmII$& $W(\7\8;\7,[\7\8];\8\8)_\rmII$&
    $W(\infty;\8,1\7;\7\8)_\rmIII$\\
    $2^{16}$&$13$&    $W(\7\8\7;;\7\7)_\rmI$& $W(\8\7\7;[\7\7];\8\7)_\rmI$&
    $W(\7\7\7;[\7\7],[\8\7];\7\8)_\rmI$\\
    &&&$W(\8\8;\8,[\7\8];\8\8)_\rmII$& $W(\7\8;\7,[\7\8],[\8\8];\7\7\7)_\rmII$\\
    $2^{17}$&$11$&    $W(\8\8\7;;\7\7)_\rmI$& $W(\7\8\7;[\7\7];\8\7)_\rmI$&
    $W(\8\7\7;[\8\7];\7\8)_\rmI$\\
    && $W(\7\7\7;[\7\7],[\7\8];\8\8)_\rmI$&& $W(\8\8;\8,[\8\8];\7\7\7)_\rmII$\\
    $2^{18}$&$19$&    $W(\7\7\8;;\7\7)_\rmI$& $W(\8\8\7;[\7\7];\8\7)_\rmI$&
    $W(\7\8\7;[\7\7],[\8\7];\7\8)_\rmI$\\
    && $W(\8\7\7;[\8\7],[\7\8];\8\8)_\rmI$&
    \multicolumn{2}{l}{$W(\7\7\7;[\7\7],[\7\8],[\8\8];\7\7\7)_\rmI$}\\
    \hline
    \end{tabular}
  \end{center}
  \caption{Normal subgroups of index up to $2^{18}$ in $\Gg$,
    contained in $K$}
  \label{table:ns}
\end{table}

\begin{proof}[Proof of Theorem~\ref{thm:normal:Gg}]
  The first two assertions are checked directly as follows. Let
  $\mathcal F$ be the set of finite-index subgroups of $\Gg$ not in
  $K$. Consider the finite quotient $Q=\Gg/\stab_6(\Gg)$, and the
  preimage $P$ of $\Gg$ defined as
  \[P = \big\langle a,b,c,d\big|\,a^2,b^2,c^2,d^2,bcd,
  \sigma^i(ad)^4,\sigma^i(adacac)^4\quad(i=0\dots5)\big\rangle.\]
  Clearly the image of $\mathcal F$ in $Q$ is at most as large as
  $\mathcal F$, and the preimage of $\mathcal F$ in $P$ is at least as
  large as $\mathcal F$. Now we use the algorithms in
  \textsc{Gap}~\cite{gap4:manual} computing the top of the lattice of
  normal subgroups for finite groups ($Q$) and finitely presented
  groups ($P$). The number of subgroups not contained in $K$ agree in
  $P$ and $Q$, so give the structure of the lattice not below $K$ in
  $\Gg$.
  
  Let now $N$ be a normal subgroup of $\Gg$, contained in $K$. If $N$
  is non-trivial, then it has finite
  index~\cite{grigorchuk:jibg}*{Corollary to Proposition~9}. It is
  easy to see that $N$ contains $C(x^2)$ and $D(x)$ for some words
  $C,D$, using for instance the congruence
  property~\cite{grigorchuk:jibg}*{Proposition~10}; therefore the
  generators of $N$ may be chosen as
  \[\{A_1(x)\cdots A_n(x)B_1(x^2)\cdots B_m(x^2),
  A'_1(x)\cdots A'_{n'}(x)B'_1(x^2)\cdots B'_{m'}(x^2),\dots,C(x^2),D(x)\},\]
  with $A_i^{(j)}<D$ and $B_i^{(j)}<C$ for all $i,j$.
  
  Taking the commutators of these generators with the appropriate
  generator among $\{a,b,c,d\}$, we shift the ranks of the $A$-terms
  up by $1$, and multiplying a generator by another we may get rid of
  all generators except $C(x^2)$ and the one with $A_1$ of smallest
  rank.
  
  We therefore consider all subgroups $W(A;B_1,\dots,B_m;C)$, and seek
  conditions on $A$, $\{B_i\}$ and $C$ so that to each normal subgroup
  in $K$ there corresponds a unique expression of the form
  $W(A;B_1,\dots,B_m;C)$.

  Let first $C$ be minimal such that $C(x^2)\in N$; then take $A$
  minimal such that for some $B_1<\dots<B_m<C$ we have
  $A(x)B_1(x^2)\cdots B_m(x^2)\in N$ . Take also $B'_1$ minimal such
  that $B'_1(x^2)\cdots B'_{m'}(x^2)\in N$ for some $B'_i$.
  
  Define the functions $M,S:\mathcal W\times2^{\mathcal
    W}\times\mathcal W\to\mathcal W$ as follows ($M$ stands for
  ``monomial'' and $S$ stands for ``squares''): Consider
  $A(x)B_1(x^2)\dots B_m(x^2)$ as an element of $\Lie_{\F[2]}(\Gg)$,
  truncated at degree $C$.  Successive commutations with generators
  $s\in\{a,b,c,d\}$, \textbf{according the the rules of
    Lemma~\ref{lem:commGg}}, give rise to other elements of
  $\Lie_{\F[2]}(\Gg)$. We stress that we use the complete computations
  of commutators, and not just those in the filtered Lie algebra.
  Define $M(A,\{B_i\})$ as the minimal word $D$ such that $D(x^2)$
  that arises in this process; if no such word occurs,
  $M(A,\{B_i\},C)=C$. Define $S(A,\{B_i\})$ as the minimal $B'_{m'}$
  such that $B'_1(x^2)\cdots B'_{m'}(x^2)$ occurs in this process; if
  no such product occurs, $S(A,\{B_i\},C)=C-1$.
  
  Now, since $M(A,\{B_i\},C)(x^2)\in N$, we necessarily have $C\le
  M(A,\{B_i\})$. Also, all $B_i$ of degree at least $B'_{m'}$ can be
  replaced by terms of lower degree $B'_1,\dots,B'_{m-1}$. This proves
  the claimed inequalities. Conversely, if there existed another
  description $A(x)\tilde B_1(x^2)\dots\tilde B_m(x^2)\in N$ for
  another choice of $\tilde B$'s, then by dividing we would obtain a
  product of $B_i(x^2)$ in $N$, contradicting $B_m<S(A,\{B_i\},C)$.
  The data $(A;B_1,\dots,B_m;C)$ subjected to the Theorem's
  constraints therefore bijectively correspond to $N$'s.

  The index of $N$ can be computed in $\Lie_{\F[2]}(\Gg)$. Seeing
  elements of $N$ as inside $\Lie$, a vector-space complement of $N$
  is spanned by all $\tilde A(x)$ of rank less than $A$, and all
  $\tilde B(x^2)$ of rank less than $S(A,\{B_i\},C)$.

  We consider finally three cases: first assume $C\le\7^{|A|}$ and
  $|B_1|\ge|A|-1$. Then $C(x^2)$ gives $\7^{|C|+1}(x^2)\7^{|C|+2}(x)$ by
  commutation with $\sigma^{|A|}(d)$, which itself gives
  $\7^{|C|}\8\7(x)$ by commutation with $a$, so we may suppose
  $A\le\7^{|C|}1\7$. Various $B_i$'s can be added, giving the
  description $(\rmI)$.
  
  Now assume $C>\7^{|A|}$. Then since $A(x)$ would produce
  $\7^{|A|}(x^2)$ by commutation with an appropriate conjugate of
  $\sigma^{|A|}(b)$, we must have $A=B_1\8$ so that the same
  commutation vanishes, giving the description $(\rmII)$.

  Finally assume we have $C\le\7^{|A|}$ and $|B_1|<|A|-1$. Then
  necessarily $A=\7^n$; and taking appropriate commutations we see that
  the normal subgroup under consideration contains
  $\7^n(x)\7^{n-1}(x^2)$. We may then replace the generator
  $A(x)B_1(x^2)\dots B_m(x^2)$ by $\7^{n-1}(x^2)B_1(x^2)\dots
  B_m(x^2)$, and obtain the description $(\rmIII)$.
\end{proof}

\begin{cor}
  Let $N$ be a normal subgroup of $\Gg$. Then $N/[N,\Gg]$ is an
  elementary $2$-group of rank $1$ or $2$, unless it is $N=\Gg$ (of
  rank $3$).
\end{cor}

\begin{cor}
  Every normal subgroup of $G$ is characteristic.
\end{cor}
\begin{proof}
  The automorphism group of $\Gg$ is determined
  in~\cite{bartholdi-s:at}: it also acts on the binary tree, and is
  \[\aut\Gg = \langle G, \8^j\7[a,d]\text{ for all }j\in\N\rangle.\]
  It then follows that $[K,\aut\Gg]=\langle\7(x),x^2\rangle^\Gg$ is a
  strict subgroup of $K$; and hence $[N,\aut\Gg]<N$ for any normal
  subgroup that is generated by expressions in $W(x)$ and $W(x^2)$ for
  words $W\in\{\7,\8\}^*$. The theorem asserts that all normal
  subgroups of $\Gg$ below $K$ have this form; it then suffices to
  check, for instance using the algorithms in
  \textsc{Gap}~\cite{gap4:manual}, that the finitely many normal
  subgroups of $\Gg$ not in $K$ are characteristic.
\end{proof}

\begin{cor}\label{cor:asympns}
  The number $b_n$ of normal subgroups of $\Gg$ of index $2^n$ starts
  as follows, and is asymptotically $n^{\log_2(3)}$. More precisely,
  we have $\liminf b_n/n^{\log_2(3)}=5^{-\log_2(3)}\approx0.078$ and
  $\limsup b_n/n^{\log_2(3)}=\frac29\approx0.222$.
  \[\begin{array}{c|cccccccccccc}
    \text{index }2^n & 2^0 & 2^1 & 2^2 & 2^3 & 2^4 & 2^5 & 2^6 & 2^7 &
    2^8 & 2^9 & 2^{10} & 2^{11}\\
    |\{N\triangleleft\Gg\}| &
    1 & 7 & 7 & 7 & 5 & 3 & 3 & 3 & 5 & 5 & 7 & 5\\ \hline
    & \rule{0em}{3ex} 2^{12} & 2^{13} & 2^{14} & 2^{15} & 2^{16} &
    2^{17} & 2^{18} & 2^{19} & 2^{20} & 2^{21} & 2^{22} & 2^{23}\\
    & 7 & 7 & 13 & 9 & 13 & 11 & 19 & 11 & 13 & 11 & 19 & 15\\ \hline
    & \rule{0em}{3ex} 2^{24} & 2^{25} & 2^{26} & 2^{27} & 2^{28} &
    2^{29} & 2^{30} & 2^{31} & 2^{32} & 2^{33} & 2^{34}\\
    & 25 & 21 & 37 & 23 & 31 & 23 & 37 & 25 & 37 & 31 & 55
  \end{array}\]
\end{cor}
\begin{proof}
  The number of subgroups of index $2^n$ behaves in a somewhat erratic
  way, but is greater when $n$ is of the form $2^k+2$, so that there
  is a maximal number of choices for $A$ and $C$, and is smaller when
  $n$ is of the form $5\cdot2^k+1$. We compute the numbers $F_k$ and
  $f_k$ of normal subgroups of $\Gg$ contained in $K$ of index $2^n$,
  with respectively $n=2^k+2$ and $n=5\cdot2^k+1$, yielding the upper
  and lower bounds. The computations are simplified by the fact that
  for these two values of $n$ there are only subgroups of type $\rmI$.
  
  Let us start by the upper bound, when $n=2^k+2$. First, for $k=2$,
  the subgroups of index $2^n$ are $W(\7;;\7)$, $W(\7;\lambda;\7)$ and
  $W(\8;;\lambda)$, giving $F_2=3$.  Then, for $k>2$, the subgroups can
  of index $2^n$ can be described as follows:
  \begin{enumerate}
  \item $W(A\8\7;\mathcal B\7;C\7)$ for all $W(A\7;\mathcal B;C)$ 
    counted in $F_{k-1}$, except when $C=\7^{k-3}$, when no subgroup
    appears in $F_k$, and when $C=\7^{k-2}$, when $C\7$ should be
    replaced by $\7^{k-3}\8$;
  \item $W(A\7;\mathcal B\8;C\8)$ for all $W(A;\mathcal B;C)$ counted
    in $F_{k-1}$, except when $C=\7^{k-3}$, when no subgroup
    appears in $F_k$, and when $C=\7^{k-2}$, when $C\8$ should be
    replaced by $\7^{k-1}$;
  \item $W(A\7;\{A\}\cup\mathcal B\8;C\8)$, with the same
    qualifications as above;
  \item $W(\7^{k-2}\8;;\7^{k-2})$.
  \end{enumerate}
  It then follows that $F_k=3(F_{k-1}-1)+1$, so $F_k=\frac293^k+1$ for
  all $k\ge2$.
  
  For the lower bound, we have $f_0=F_2=3$; and for $k>0$, when
  $n=5\cdot2^k+1$, the subgroups can of index $2^n$ can be described
  as follows:
  \begin{enumerate}
  \item $W(A\8\8;\mathcal B\7;C\7)$ for all $W(A\8;\mathcal B;C)$ 
    counted in $f_{k-1}$;
  \item $W(A\7\8;\mathcal B\8;C\8)$ for all $W(A\8;\mathcal B;C)$
    counted in $f_{k-1}$;
  \item $W(A\7\8;\{A\7\}\cup\mathcal B\8;C\8)$, with the same
    qualifications as above;
  \item $W(\8^k\7;;\7^{k+1})$ and $W(\8^k\7;\8^k;\7^{k+1})$.
  \end{enumerate}
  It then follows that $f_k=3(f_{k-1}-2)+2$, so $F_k=3^k+2$ for all
  $k\ge0$.
  
  In summary, the number of normal subgroups of index $2^n$ oscillates
  between $3^{\log_2(\frac{n-1}{5})}+2$ and $\frac293^{\log_2(n-2)}+1$
  for $n\ge6$ (when all normal subgroups of $\Gg$ are contained in
  $K$). These bounds give respectively
  $5^{-\log_2(3)}(n-1)^{\log_2(3)}$ and $\frac29(n-2)^{\log_2(3)}$.
\end{proof}

Note also the following curiosity:
\begin{cor}\label{cor:allodd}
  The number of normal subgroups of index $r$ of $G$ is odd for all
  $r$'s a power of $2$, and even (in fact, $0$) for all other $r$.
\end{cor}
\noindent(The same congruence phenomenon holds for the group
$C_2*C_3$, as observed by Thomas M\"uller~\cite{muller:subgroups})
\begin{proof}
  The proof follows from the description of
  Theorem~\ref{thm:normal:Gg}.  Assume $r=2^k$. To determine the
  parity of the number of subgroups of index $r$, it suffices to
  consider which $W(A;\mathcal B;C)$ expressions have no choices for
  $\mathcal B$. These are precisely the $W(A;;\7^n)_\rmI$ with
  $2^{n+1}<\#A\le5\cdot2^n$, the $W(\7^n\8\7;;C)_\rmI$ with
  $2^n<\#C\le2^{n+1}$ and the $W(\infty;\8^n,C-1;C)_\rmIII$ with
  $2^{n+1}+1<\#C\le3\cdot2^n+1$.

  Now these last two families yield a subgroup for precisely the same
  values of $k$, namely those satisfying $6\cdot2^j+2\le
  k\le7\cdot2^j+1$, and therefore contribute nothing modulo $2$. The
  first family contributes a subgroup for all $k$.
\end{proof}


\subsection{Normal subgroups in $\GS$} The normal subgroup growth
of $\GS$ is much larger. As a crude lower bound, consider the
quotient $A=\gamma_k(\GS)/\gamma_{k+1}(\GS)$ for
$k=\frac12(\alpha_{2n+1}+1)$. It is abelian of rank $2^n$, and the
index of $\gamma_k(\GS)$, respectively $\gamma_{k+1}(\GS)$, is
$3^{3^{2n-1}\pm2^{n-1}+1}$.

In the vector space $\F[3]^j$, there are roughly $3^{\binom j2}$
subspaces; so $A$ has about $3^{4^n}$ subgroups
$S=N/\gamma_{k+1}(\GS)$, each of them giving rise to a subgroup $N$ of
index roughly $3^{9^n}$.

It then follows that the number of normal subgroups of $\GS$ of
index $3^n$ is at least $3^{n^{\log_3(2)}}$, a function intermediate
between polynomial and exponential growth. More precise estimations of
the normal subgroup growth of $\GS$ will be the topic of a future
paper.

\subsection*{Acknowledgments}
I wish to express my immense gratitude to the referee who helped me
clarify many parts of the present and forthcoming paper.

\end{document}